# UPPER BOUNDS FOR SPATIAL POINT PROCESS APPROXIMATIONS

By Dominic Schuhmacher[1]

*University of Zürich*

We consider the behavior of spatial point processes when subjected to a class of linear transformations indexed by a variable $T$. It was shown in Ellis [*Adv. in Appl. Probab.* **18** (1986) 646–659] that, under mild assumptions, the transformed processes behave approximately like Poisson processes for large $T$. In this article, under very similar assumptions, explicit upper bounds are given for the $d_2$-distance between the corresponding point process distributions. A number of related results, and applications to kernel density estimation and long range dependence testing are also presented. The main results are proved by applying a generalized Stein–Chen method to discretized versions of the point processes.

**1. Introduction.** Let $D_1, D_2 \in \mathbb{N} = \{1, 2, 3, \ldots\}$ and $D = D_1 + D_2$. Consider a point process $\xi$ on $\mathbb{R}^D = \mathbb{R}^{D_1} \times \mathbb{R}^{D_2}$, which has expectation measure $\nu$ and meets three conditions, namely, absolute continuity of $\nu$ with a mild restriction on the density, an orderliness condition in the $\mathbb{R}^{D_1}$-directions and a mixing condition in the $\mathbb{R}^{D_2}$-directions (formal versions of these conditions can be found at the end of this section). Let $\eta$ be a Poisson process with the same expectation measure and let $\theta_T : \mathbb{R}^D \to \mathbb{R}^D$ be the linear transformation that stretches the first $D_1$ coordinates by a factor $w(T)^{1/D_1}$ and compresses the last $D_2$ coordinates by a factor $T^{1/D_2}$, that is, for $T \in \mathbb{R}$, $T \geq 1$, we set

$$\theta_T(\mathbf{s}, \mathbf{t}) := \left(w(T)^{1/D_1}\mathbf{s}, \frac{1}{T^{1/D_2}}\mathbf{t}\right) \quad \text{for all } (\mathbf{s}, \mathbf{t}) \in \mathbb{R}^{D_1} \times \mathbb{R}^{D_2} = \mathbb{R}^D,$$

where $w(T) \to \infty$ and $w(T) = O(T)$ for $T \to \infty$. In particular, we usually write $\tilde{\theta}_T$ instead of $\theta_T$ if our stretch factor is $T^{1/D_1}$.

Received August 2001; revised October 2003.
[1]Supported in part by Schweizerischer Nationalfonds, Projekt Nr. 20-61753.00.

Most of the time we will restrict our transformed processes $\xi\theta_T^{-1}$ and $\eta\theta_T^{-1}$ to a bounded cube $J := [-1, 1)^D$ and denote by $J_T := \theta_T^{-1}(J)$ the pre-image of $J$, but sometimes the bigger cuboids $\tilde{J}_T := \tilde{\theta}_T(J_T) = [-(\frac{T}{w(T)})^{1/D_1}, (\frac{T}{w(T)})^{1/D_1})^{D_1} \times [-1, 1)^{D_2}$ instead of $J$ are more useful.

A consequence of what Ellis (1986) showed is that, for bounded measurable functions $f_T : J \to \mathbb{R}$ with $\|f_T\|_\infty = O(\sqrt{w(T)/T})$, the distributions of $\int_J f_T \, d(\xi\theta_T^{-1})$ and $\int_J f_T \, d(\eta\theta_T^{-1})$ get more and more alike as $T \to \infty$; or, more precisely, that the difference between their characteristic functions converges uniformly to zero on every compact subset of $\mathbb{R}$ as $T \to \infty$. Therefore, there is hope that $d(\mathcal{L}(\xi\theta_T^{-1}|_J), \mathcal{L}(\eta\theta_T^{-1}|_J))$ can be shown to be small for large $T$ if we choose for $d$ a probability distance between distributions of point processes which metrizes a topology that is equal to or not too much finer than the weak topology (i.e., the topology of convergence in distribution).

Our choice for $d$ will be the $d_2$-distance [see Barbour, Holst and Janson (1992), Section 10.2], which, besides meeting the aforementioned requirement, has a number of other useful properties; it is rather easy to handle, and bounds on $d_2(\mathcal{L}(\xi_1), \mathcal{L}(\xi_2))$ for point processes $\xi_1, \xi_2$ imply bounds on $|\mathbb{E}f(\xi_1) - \mathbb{E}f(\xi_2)|$ for a number of desirable functions $f$. The $d_2$-distance can be constructed as two Wasserstein distances, one on top of the other, in the following way. Consider a compact set $\mathcal{X} \subset \mathbb{R}^D$ and write $\mathcal{M}_p$ for the space of point measures on $\mathcal{X}$. Let $d_0$ be the usual Euclidean distance on $\mathbb{R}^D$, but bounded by 1, and $\mathcal{F}_1 := \{k : \mathcal{X} \to \mathbb{R}; |k(x_1) - k(x_2)| \leq d_0(x_1, x_2)\}$. Define the $d_1$-*distance* (w.r.t. $d_0$) between point measures $\rho_1, \rho_2 \in \mathcal{M}_p$ by

$$d_1(\rho_1, \rho_2) := \begin{cases} 1, & \text{if } |\rho_1| \neq |\rho_2|, \\ \frac{1}{|\rho_1|} \sup_{k \in \mathcal{F}_1} \left| \int k \, d\rho_1 - \int k \, d\rho_2 \right|, & \text{if } |\rho_1| = |\rho_2| \geq 1, \\ 0, & \text{if } |\rho_1| = |\rho_2| = 0, \end{cases}$$

where $|\rho_i| := \rho_i(\mathcal{X}) < \infty$. It can be seen that $(\mathcal{M}_p, d_1)$ is a complete, separable metric space and that $d_1$ is bounded by 1. Furthermore, the Kantorovich–Rubinstein theorem [see Dudley (1989), Section 11.8] when $|\rho_1| = |\rho_2| =: n \geq 1$ yields that

$$(1.1) \qquad d_1(\rho_1, \rho_2) = \min_{\pi \in S_n} \left[ \frac{1}{n} \sum_{i=1}^n d_0(x_{1,i}, x_{2,\pi(i)}) \right],$$

where $S_n$ is the set of permutations of $\{1, 2, \ldots, n\}$. Now let $\mathcal{F}_2 := \{f : \mathcal{M}_p \to \mathbb{R}; |f(\rho_1) - f(\rho_2)| \leq d_1(\rho_1, \rho_2)\}$ and define the $d_2$-*distance* (w.r.t. $d_0$) between probability measures $P$ and $Q$ on $\mathcal{M}_p$ (distributions of point processes on $\mathcal{X}$) by

$$d_2(P, Q) := \sup_{f \in \mathcal{F}_2} \left| \int f \, dP - \int f \, dQ \right|.$$

By the Kantorovich–Rubinstein theorem, one obtains that

$$(1.2) \qquad d_2(P,Q) = \min_{\substack{\xi_1 \sim P \\ \xi_2 \sim Q}} \mathbb{E} d_1(\xi_1, \xi_2)$$

[the minimum is attained, because $(\mathcal{M}_p, d_1)$ is complete, see Rachev (1984)]. Furthermore, because of the bound on the $d_1$-distance, the $d_2$-distance can also be interpreted as a variant of a bounded Wasserstein distance (see below). Hence, Theorem 11.3.3 in Dudley (1989) yields that $d_2$ metrizes the weak convergence of point process distributions; or, in other words, for point processes $\xi, \xi_1, \xi_2, \ldots$ on $\mathcal{X}$, we have

$$(1.3) \qquad \xi_n \xrightarrow{\mathcal{D}} \xi \quad \text{iff } d_2(\mathcal{L}(\xi_n), \mathcal{L}(\xi)) \to 0,$$

where the convergence in distribution for point processes is defined in the usual sense [see Kallenberg (1986), Section 4.1]. The fact that is crucial here is that, for $d_0$ as defined, the topology generated by the metric $d_1$ on $\mathcal{M}_p$ is equal to the vague topology, which is used for the definition of convergence in distribution for point processes.

$d_2$ is the distance that we are mainly interested in, but we will also deal with two other probability distances; namely, on the one hand, the *total variation distance* between distributions $\mu_1$ and $\mu_2$ on $\mathbb{Z}_+$, which is defined as

$$d_{\mathrm{TV}}(\mu_1, \mu_2) := \sup_{A \subset \mathbb{Z}_+} |\mu_1(A) - \mu_2(A)|$$

and can be equivalently written in the form

$$(1.4) \qquad d_{\mathrm{TV}}(\mu_1, \mu_2) = \min_{\substack{X_1 \sim \mu_1 \\ X_2 \sim \mu_2}} \mathbb{P}[X_1 \neq X_2];$$

and, on the other hand, the *bounded Wasserstein distance* between distributions $\tilde{\mu}_1$ and $\tilde{\mu}_2$ on $\mathbb{R}$, which is defined as

$$d_{\mathrm{BW}}(\tilde{\mu}_1, \tilde{\mu}_2) := \sup_{f \in \mathcal{F}_{\mathrm{BW}}} \left| \int_{\mathbb{R}} f \, d\tilde{\mu}_1 - \int_{\mathbb{R}} f \, d\tilde{\mu}_2 \right|,$$

where

$$\mathcal{F}_{\mathrm{BW}} := \{f : \mathbb{R} \to \mathbb{R}; |f(x) - f(y)| \leq |x - y| \text{ and } |f(x)| \leq \tfrac{1}{2} \text{ for } x, y \in \mathcal{X}\},$$

the set of Lipschitz continuous functions with constant 1 that are bounded by $\frac{1}{2}$. For equivalent expressions and properties see Barbour, Holst and Janson (1992), Appendix A.1 for the total variation distance and Dudley (1989), Section 11.3 for the bounded Wasserstein distance.

It will be the main goal of our endeavors to find upper estimates for the distance $d_2(\mathcal{L}(\xi\theta_T^{-1}|_J), \mathcal{L}(\eta\theta_T^{-1}|_J))$ (see Section 2.2), but explicit upper bounds will also be computed for $d_{\mathrm{TV}}(\mathcal{L}(\xi\theta_T^{-1}(J)), \mathcal{L}(\eta\theta_T^{-1}(J)))$ (Section 2.3), $d_2(\mathcal{L}(\xi\tilde{\theta}_T^{-1}|_{\tilde{J}_T}), \mathcal{L}(\eta\tilde{\theta}_T^{-1}|_{\tilde{J}_T}))$ (Section 2.4) and $d_2(\mathcal{L}(\xi\tilde{\theta}_T^{-1}|_{\tilde{J}_T}), \mathrm{Po}(\nu'|_{\tilde{J}_T}))$ for an appropriate $T$-independent measure $\nu'$ on $\mathbb{R}^D$ (Section 2.5). Throughout the article we use $\mathrm{Po}(\nu')$ to denote the Poisson distribution with parameter $\nu'$ if $\nu'$ is a positive real number and to denote the distribution of the Poisson process with parameter measure $\nu'$ if $\nu'$ is a boundedly finite measure.

In Section 3 we present some applications of our results. Most importantly, we calculate an upper bound for the bounded Wasserstein distance between the distribution of a kernel estimate of the density of $\nu$ at a certain point and the actual value of the density at that point. Furthermore, we briefly describe an application to testing for long range dependence.

Apart from the paper of Ellis (1986), which provided the initial motivation for many of the theorems in this article, stretched point processes have also been investigated in the context of light traffic analysis for queues and in other, similar topics: see, for example, Borovkov (1996) and the references therein. These authors, however, were interested in the quite different question of finding asymptotic expansions for the expectation of functionals of purely stretched marked point processes, which vanish in the limit on every compact set; our procedure, in contrast, leads to point processes with, essentially, a stable or increasing number of points in every compact set.

We conclude this section by having a detailed look at the three conditions for the point process $\xi$.

CONDITION 1 (Absolute continuity of the expectation measure). Let $\mu = \mu_1 \otimes \mu_2$, where $\mu_1 := \lambda^{D_1}$ is the Lebesgue measure on $\mathbb{R}^{D_1}$, and either $\mu_2 := \lambda^{D_2}$ is the Lebesgue measure on $\mathbb{R}^{D_2}$ or $\mu_2 := \mathcal{H}_0^{D_2}$ is the counting measure on $\mathbb{Z}^{D_2} + \frac{1}{2}\mathbf{1} \subset \mathbb{R}^{D_2}$.

Then we require that $\nu \ll \mu$ with a Radon–Nikodym density $p$, such that $\kappa \in \mathbb{R}_+$ exists with

$$\kappa_T := \sup_{(\mathbf{s},\mathbf{t}) \in J_T} p(\mathbf{s},\mathbf{t}) \leq \kappa \qquad \text{for all } T \geq 1.$$

In the same way, we choose $\iota \in \mathbb{R}_+$ with

$$\iota_T := \inf_{(\mathbf{s},\mathbf{t}) \in J_T} p(\mathbf{s},\mathbf{t}) \geq \iota \qquad \text{for all } T \geq 1.$$

(For the asymptotic result it is enough, of course, to assume both statements only for all $T$ bigger than some $T_0 \geq 1$.)

CONDITION 2 (Orderliness). There is a continuous function $\breve{\alpha}\colon\mathbb{R}_+ \to \mathbb{R}_+$ with $\breve{\alpha}(0) = 0$, such that for every rectangle $C := [\mathbf{a}, \mathbf{b}) \times [\mathbf{c}, \mathbf{d})$ with $\mathbf{a}, \mathbf{b} \in \mathbb{R}^{D_1}$, $\mathbf{a} \leq \mathbf{b}$, and $\mathbf{c}, \mathbf{d} \in \mathbb{R}^{D_2}$, $\mathbf{c} \leq \mathbf{d}$, we have

$$\mathbb{E}[(\xi(C))^2 \mathbb{1}_{\{\xi(C) \geq 2\}}] \leq v\breve{\alpha}(v),$$

where

$$v := v(C) = \mu_1([\mathbf{a}, \mathbf{b}))\mu_2([\mathbf{c}, \mathbf{d} + \mathbf{1})).$$

For the third condition, there are different versions that can be considered. According to the type of mixing we are interested in, we write this condition as $3x$, where $x \in \{\beta, \rho, \varphi\}$:

CONDITION $3x$ ($x$-mixing property). For every interval $[\mathbf{a}, \mathbf{b}) \subset \mathbb{R}^{D_1}$, $\mathbf{a} < \mathbf{b}$, there is a decreasing function $\breve{\beta} := \breve{\beta}_{\mathbf{a},\mathbf{b}} \colon \mathbb{R}_+ \to \mathbb{R}_+$ with the two following properties:

(a) $\breve{\beta}(u) = o(\frac{1}{u^{D_2/2}})$ for $u \to \infty$.
(b) If $\mathbf{c}, \mathbf{d} \in \mathbb{R}^{D_2}$ with $\mathbf{c} < \mathbf{d}$, $t \in \mathbb{R}_+$ and the $\sigma$-fields $\mathcal{F}_{\text{int}}$ and $\mathcal{F}_{\text{ext}}$ are defined as $\mathcal{F}_{\text{int}} := \sigma(\xi|_{[\mathbf{a},\mathbf{b})\times[\mathbf{c},\mathbf{d})})$ and $\mathcal{F}_{\text{ext}} := \sigma(\xi|_{[\mathbf{a},\mathbf{b})\times[\mathbf{c}-t\mathbf{1},\mathbf{d}+t\mathbf{1})^c})$, then

$$x(\mathcal{F}_{\text{int}}, \mathcal{F}_{\text{ext}}) \leq \breve{\beta}(t),$$

where $x$ is one of the three mixing coefficients $\beta, \rho$ or $\varphi$ with

$$\beta(\mathcal{F}_{\text{int}}, \mathcal{F}_{\text{ext}}) := \mathbb{E}\operatorname*{ess\,sup}_{B \in \mathcal{F}_{\text{ext}}} |\mathbb{P}(B|\mathcal{F}_{\text{int}}) - \mathbb{P}(B)|,$$

$$\rho(\mathcal{F}_{\text{int}}, \mathcal{F}_{\text{ext}}) := \sup_{\substack{X \in L_2(\mathcal{F}_{\text{int}}) \\ Y \in L_2(\mathcal{F}_{\text{ext}})}} |\operatorname{corr}(X, Y)|,$$

$$\varphi(\mathcal{F}_{\text{int}}, \mathcal{F}_{\text{ext}}) := \sup_{\substack{A \in \mathcal{F}_{\text{int}} \\ B \in \mathcal{F}_{\text{ext}}}} |\mathbb{P}(B|A) - \mathbb{P}(B)|.$$

In the following we suppress the indication of the interval $[\mathbf{a}, \mathbf{b})$ and write simply $\breve{\beta}$. The corner points $\mathbf{a}$ and $\mathbf{b}$ are to be chosen appropriately; for example, $\mathbf{a} = -\sup_{T \geq 1}(\frac{1}{w(T)})^{1/D_1} \cdot \mathbf{1}$, $\mathbf{b} = \sup_{T \geq 1}(\frac{1}{w(T)})^{1/D_1} \cdot \mathbf{1}$ is always an appropriate choice.

No further explanation is needed for the first condition. It simply states the absolute continuity of the expectation measure with respect to what is basically Lebesgue measure, with a mild condition on the density. The fact that we admit the counting measure for the $D_2$-part of the reference measure $\mu$ allows us to apply our future estimates to (mixing) sequences of certain $\mathbb{R}^{D_1}$-valued point processes. In order to simplify certain formulas, we will always tacitly assume that $T \in \{n^{D_2}; n \in \mathbb{N}\}$ if $\mu_2$ is the counting measure.

The second condition is a form of orderliness in the $\mathbb{R}^{D_1}$-directions. For a detailed account of orderliness, see Daley (1974). For what we are interested in here, it is enough to understand that the upper bound for $\mathbb{E}[(\xi(C))^2 \mathbb{1}_{\{\xi(C) \geq 2\}}]$ implies that

$$4\mathbb{P}[\xi(C) \geq 2] \leq v\breve{\alpha}(v),$$

and that Condition 2 implies the simplicity of $\xi$ (i.e., $\mathbb{P}[\xi(\{x\}) \leq 1 \ \forall x \in \mathbb{R}^D] = 1$). The latter implication is due to Theorem 2.6 in Kallenberg (1986).

The various versions of the third condition are mixing conditions of different strength. It can be seen [Doukhan (1994)] that

$$\beta(\mathcal{B}, \mathcal{C}) \leq \varphi(\mathcal{B}, \mathcal{C}),$$
$$\rho(\mathcal{B}, \mathcal{C}) \leq 2\varphi^{1/2}(\mathcal{B}, \mathcal{C})\varphi^{1/2}(\mathcal{C}, \mathcal{B})$$

for arbitrary $\sigma$-fields $\mathcal{B}, \mathcal{C} \subset \mathcal{F}$ on some common probability space $(\Omega, \mathcal{F}, \mathbb{P})$. Thus, the concept of $\varphi$-mixing is the strongest of the three, followed by the $\beta$-mixing and $\rho$-mixing concepts, which are not generally comparable with each other, although from an empirical point of view, $\beta$-mixing often turns out to be the stronger of the two. Two mixing concepts that are not treated here are $\alpha$-mixing, which would be weaker, and $\psi$-mixing, which would be stronger than any of the three mentioned concepts [see Doukhan (1994)]. The kind of mixing used in Ellis (1986) is $\rho$-mixing. However, it is important to notice that we need a stronger mixing condition, in the sense that the set underlying the $\sigma$-field $\mathcal{F}_{\text{ext}}$ may enclose the set underlying the $\sigma$-field $\mathcal{F}_{\text{int}}$ from all of the $2D_2$ possible directions of the $\mathbb{R}^{D_2}$. As partial compensation, the order we need for the convergence of our mixing coefficient to zero is only half the order that was needed for Ellis' result, and what is more, we could actually manage with a mixing condition where the $\sigma$-fields $\mathcal{F}_{\text{ext}}$ and $\mathcal{F}_{\text{int}}$ are quite a bit smaller (namely, generated by the numbers of points of $\xi$ in the corresponding discretization cuboids that we will need for the proof).

**2. The main results.** The results given within this section have somewhat similar flavor, and their proofs all follow the same path; first discretizing the point processes and then applying a local Stein theorem. An outline of this method can be found in Section 2.1; thereafter, in Sections 2.2–2.5 the different results are presented. A detailed, self-contained proof is given only for Theorem 2.A; for the other statements the necessary adaptations are given.

2.1. *The approach.* All statements in Section 2 are about upper bounds for distances between the distribution of a transformed $\xi$-process and the distribution of a transformed Poisson process (or a function of the respective process, as in Section 2.3). For the sake of clarity of presentation, we formulate the ideas of the proof only for $d_2(\mathcal{L}(\xi\theta_T^{-1}|_J), \mathcal{L}(\eta\theta_T^{-1}|_J))$. However, except for the obvious changes in notation (like writing $\xi\tilde{\theta}_T^{-1}|_{\tilde{J}_T}$ instead of $\xi\theta_T^{-1}|_J$ in Section 2.4), the arguments presented here can be applied *literally* (or almost literally in the case of Section 2.3) to calculate the presented upper bounds for any of the distances appearing in this section.

As mentioned before, our basic strategy of proof is to discretize $\xi\theta_T^{-1}$ and $\eta\theta_T^{-1}$ (in general, the point processes involved) and then apply an estimate, obtained by a generalized version of the Stein–Chen method, to the discretized point processes (in fact, the classic Stein–Chen method will be enough for Section 2.3, where only the numbers of points are involved). The corresponding estimate can be found in the Appendix.

The discretizations are carried out as follows. For every $T \geq 1$ and for $h(T) \geq 1$, set $n_1 := \lceil h(T)^{1/D_1} \rceil - 1$ and $n_2 := \lceil T^{1/D_2} \rceil - 1$, where $\lceil x \rceil$ denotes, for any $x \in \mathbb{R}$, the smallest integer $z \geq x$. We subdivide $J_T$ into smaller "discretization cuboids" $C_{\mathbf{kl}}$ with lengths 1 in the $\mathbb{R}^{D_2}$-directions and widths $\frac{1}{(w(T)h(T))^{1/D_1}}$ in the $\mathbb{R}^{D_1}$-directions, whenever the $C_{\mathbf{kl}}$ are not too close to the boundary of $J_T$. Here $h(T)$ can be thought of as order of the number of discretization cuboids in the $\mathbb{R}^{D_1}$-directions [there are $2\lceil h(T)^{1/D_1} \rceil$ in every dimension of $\mathbb{R}^{D_1}$]. To be more precise, we set, for every $T \geq 1$,

$$C_{\mathbf{kl}} := C_{\mathbf{kl}}^{(T)}$$
$$:= \left(\prod_{r=1}^{D_1}\left[-\frac{n_1}{(w(T)h(T))^{1/D_1}} + \frac{k_r - 1}{(w(T)h(T))^{1/D_1}},\right.\right.$$
$$\left.-\frac{n_1}{(w(T)h(T))^{1/D_1}} + \frac{k_r}{(w(T)h(T))^{1/D_1}}\right)$$
$$\times \prod_{s=1}^{D_2}[-n_2 + (l_s - 1), -n_2 + l_s)\right) \cap J_T$$

for all $\mathbf{k} = (k_1, k_2, \ldots, k_{D_1}) \in \{0, 1, \ldots, 2n_1 + 1\}^{D_1}$ and $\mathbf{l} = (l_1, l_2, \ldots, l_{D_2}) \in \{0, 1, \ldots, 2n_2 + 1\}^{D_2}$, so that $J_T = \dot{\bigcup}_{\mathbf{k},\mathbf{l}} C_{\mathbf{kl}}^{(T)}$. Note that in order to reduce the complexity of presentation, we will make use of simplified notations for multi-indices that should be obvious in their meaning. For instance, we write, in short, $\sum_{\mathbf{k}=0}^{2n_1+1} a_\mathbf{k}$ instead of $\sum_{\mathbf{k}:\, k_1, \ldots k_r = 0}^{2n_1+1} a_\mathbf{k}$ or $\mathbf{k} \in \{0, 1, \ldots, 2n_1 + 1\}$ instead of $\mathbf{k} \in \{0, 1, \ldots, 2n_1 + 1\}^{D_1}$. Also, where not stated otherwise, the ranges of the indices in expressions like $\sum_{\mathbf{k},\mathbf{l}}$ or $\bigcup_{\mathbf{k},\mathbf{l}}$ are given by $\mathbf{k} \in \{0, 1, \ldots, 2n_1 + 1\}$, $\mathbf{l} \in \{0, 1, \ldots, 2n_2 + 1\}$. Some more notation is needed. We denote by $\alpha_{\mathbf{kl}}$ the centre of $C_{\mathbf{kl}}$ and define in the image space of the transformation $\theta_T$

$$R_{\mathbf{kl}} := R_{\mathbf{kl}}^{(T)} := \theta_T(C_{\mathbf{kl}}^{(T)})$$

$$= \prod_{r=1}^{D_1} \biggl[ -\frac{n_1}{h(T)^{1/D_1}} + \frac{k_r-1}{h(T)^{1/D_1}}, -\frac{n_1}{h(T)^{1/D_1}} + \frac{k_r}{h(T)^{1/D_1}} \biggr)$$

$$\times \prod_{s=1}^{D_2} \biggl[ -\frac{n_2}{T^{1/D_2}} + \frac{l_s-1}{T^{1/D_2}}, -\frac{n_2}{T^{1/D_2}} + \frac{l_s}{T^{1/D_2}} \biggr)$$

for all $\mathbf{k}, \mathbf{l}$ and write $\rho_{\mathbf{kl}}$ for the centre of $R_{\mathbf{kl}}$ [correspondingly, we use $\tilde{R}_{\mathbf{kl}} := \tilde{\theta}_T(C_{\mathbf{kl}}^{(T)})$ and $\tilde{\rho}_{\mathbf{kl}}$ in Section 2.4].

The discretization $\Xi$ of the point process $\xi$ is obtained by setting a point in the middle of every discretization cuboid $C_{\mathbf{kl}}$ which contains any points of $\xi$. Formally, we set

$$I_{\mathbf{kl}} := I_{\mathbf{kl}}^{(T)} := \mathbb{1}_{\{\xi(C_{\mathbf{kl}}) \geq 1\}}, \qquad p_{\mathbf{kl}} := \mathbb{E} I_{\mathbf{kl}} \qquad \text{for all } \mathbf{k}, \mathbf{l},$$

$$W := W^{(T)} := \sum_{\mathbf{k},\mathbf{l}} I_{\mathbf{kl}}, \qquad \lambda := \mathbb{E} W = \sum_{\mathbf{k},\mathbf{l}} p_{\mathbf{kl}},$$

and define $\Xi$ as

$$\Xi := \sum_{\mathbf{k},\mathbf{l}} I_{\mathbf{kl}} \delta_{\alpha_{\mathbf{kl}}}.$$

The error we make in the transition from $\xi \theta_T^{-1}|_J$ to $\Xi \theta_T^{-1}$ in terms of the $d_2$-distance (with a slight alteration, the argument holds also for the $d_{\mathrm{TV}}$-distance between the numbers of points; see Section 2.3) is small for large $T$, because, on the one hand, the orderliness condition (Condition 2) takes care that the probability of two points within the same discretization cuboid (and, as a consequence, of any point vanishing in the transition) is small, and, on the other hand, we have chosen our discretization in such a way that we only have to move points by a $d_0$-distance of, at most, half a body diagonal of a discretization cuboid $R_{\mathbf{kl}}$ ($\tilde{R}_{\mathbf{kl}}$ in Section 2.4) in the image space, which is small for large $T$ as well.

As a discretization (at least "in distribution") of the Poisson point process $\eta$, we take

$$\mathrm{H} := \sum_{\mathbf{k},\mathbf{l}} U_{\mathbf{kl}} \delta_{\alpha_{\mathbf{kl}}},$$

where $U_{\mathbf{kl}}$ are arbitrary independent $\mathrm{Po}(p_{\mathbf{kl}})$-distributed random variables for $0 \leq \mathbf{k} \leq 2n_1 + 1$, $0 \leq \mathbf{l} \leq 2n_2 + 1$. Again, the error we make in the transition from $\eta \theta_T^{-1}|_J$ to $\mathrm{H} \theta_T^{-1}$ is small for reasons quite similar to those stated above for the transition from $\xi \theta_T^{-1}|_J$ to $\Xi \theta_T^{-1}$ (note that the two discretizations were not realized in the same way, and that we have to argue a little more carefully in Section 2.5, where a limiting Poisson process that is independent of $T$ is considered).

We then have an indicator point process $\Xi$ with a local dependence property (stemming from the mixing Condition $3x$) and a discrete Poisson point process with the appropriate intensity measure, so that we are in the position to apply the local Stein Theorem A.D for point processes (or, in case of Section 2.3, Theorem A.A for sums of indicators), which in each case yields the stated result.

There is one point about the refinement of our discretization that is worth noting. In our main $\rho$-mixing case we retain the highest possible flexibility by introducing the variable $h(T)$. Although it will often turn out to be a natural and relatively good choice to set $h(T) := T$, doing so is, in many cases, not optimal. The optimal choice of $h(T)$ depends on the specific orderliness and mixing conditions that can be obtained for $\xi$. The weaker the orderliness condition [the slower $\breve{\alpha}(v)$ goes to zero for $v \to 0$], the higher the optimal $h(T)$ will be; conversely (and somewhat surprisingly at the moment), the weaker the mixing condition [the slower $\breve{\beta}(u)$ goes to zero for $u \to \infty$], the lower the optimal $h(T)$ will be. In contrast, no such considerations are necessary for the discretization in the $\mathbb{R}^{D_2}$-directions. A discretization cuboid length of 1 can easily be seen to be both natural and optimal. A length of higher order in $T$ only increases the distance, by which we have to move points for discretizing, a length of lower order in $T$ increases the number of discretization cuboids without changing the order of the length that the orderliness condition "sees" [i.e., without changing $v(C_{\mathbf{kl}})$ with $v$ as in Condition 2].

2.2. *The $d_2$-distance between the point processes.* In this section the $d_2$-distance between the transformed point processes $\xi \theta_T^{-1}|_J$ and $\eta \theta_T^{-1}|_J$ is considered. In all the results we use the notation $O(f_1(T), \ldots, f_j(T))$ as short hand for $O(\max\{f_1(T), \ldots, f_j(T)\})$.

2.2.1. *Results.*

THEOREM 2.A ("The principal theorem"). *Suppose that the prerequisites of Section 1 hold, including the Conditions 1, 2 and $3\rho$, and let $\iota > 0$.*

*Then we obtain for arbitrary $m := m(T) \in \mathbb{Z}_+$ and $h(T) \geq 1$ for every $T \geq 1$:*

$$d_2(\mathcal{L}(\xi\theta_T^{-1}|_J), \mathcal{L}(\eta\theta_T^{-1}|_J))$$
$$= O\bigg(\frac{1}{h(T)^{1/D_1}}, \frac{1}{T^{1/D_2}}, \log^\uparrow\bigg(\frac{T}{w(T)}\bigg)\frac{m^{D_2}+1}{w(T)}, \frac{T}{w(T)}\breve{\alpha}\bigg(\frac{2^{D_2}}{w(T)h(T)}\bigg),$$
$$\log^\uparrow\bigg(\frac{T}{w(T)}\bigg)\breve{\alpha}\bigg(\frac{2^D(2m+1)^{D_2}}{w(T)}\bigg), \sqrt{Th(T)}\breve{\beta}(m)\bigg)$$
$$\text{for } T \to \infty,$$

*where we write* $\log^{\uparrow}(x) := 1 + (\log(x) \vee 0)$ *for* $x > 0$.

For a quantitative form of the upper bound see (2.10) and (2.11) at the end of the proof. Note that the powers of 2 and 5 that appear in these inequalities have been chosen (for the convenience of calculations) to be unnecessarily large and might be dramatically improved.

One now might ask the question under what conditions the $d_2$-distance converges to zero.

COROLLARY 2.B (Convergence to zero in Theorem 2.A). *Suppose that the prerequisites of Theorem 2.A hold. Furthermore, suppose that* $w(T) \geq kT^{\delta}$ *for* $k > 0$, $\delta \in (0,1]$ *and that*

$$\breve{\alpha}(v) = O(v^r) \qquad \text{for } v \to 0 \text{ with } r > 0,$$

$$\breve{\beta}(u) = O\left(\frac{1}{u^{(1+s)D_2/2}}\right) \qquad \text{for } u \to \infty \text{ with } 1+s > \max\left(\frac{1-\delta}{\delta}\frac{1+r}{r}, \frac{1}{\delta}\right).$$

*Then*

$$d_2(\mathcal{L}(\xi \theta_T^{-1}|_J), \mathcal{L}(\eta \theta_T^{-1}|_J)) \to 0 \quad \text{for } T \to \infty.$$

REMARK 2.C (Convergence to zero, simplified).

(a) By adjusting $m$ and $h(T)$ to the function $\breve{\beta}$ it can be shown easily that for $w(T) \asymp T$, the convergence $d_2(\mathcal{L}(\xi \theta_T^{-1}|_J), \mathcal{L}(\eta \theta_T^{-1}|_J)) \to 0$ holds under the general prerequisits of Theorem 2.A. This is consistent with Corollary 2.B for $\delta = 1$ (note that the requirements for the functions $\breve{\alpha}$ and $\breve{\beta}$ are a bit stronger in Corollary 2.B).

(b) From Corollary 2.B follows that for arbitrary $\delta \in (0,1]$ and for $r > \frac{1-\delta}{1+\delta}$, $1+s > \frac{2}{\delta}$, we have $d_2(\mathcal{L}(\xi \theta_T^{-1}|_J), \mathcal{L}(\eta \theta_T^{-1}|_J)) \to 0$ for $T \to \infty$. These simpler, but stronger requirements on the functions $\breve{\alpha}$ and $\breve{\beta}$ reflect the case where we refrain from adapting $h(T)$ to the concrete problem and simply set $h(T) = T$.

In the principal Theorem 2.A, it may seem a little unsatisfactory that our "discretization depth" $h(T)$ in the $\mathbb{R}^{D_1}$-directions appears in the term $\sqrt{Th(T)}\breve{\beta}(m)$, which stems from the mixing condition in the $\mathbb{R}^{D_2}$-directions, and that, in fact, a finer discretization could increase the overall upper bound we get for the $d_2$-distance. Whereas it might well be that the factor $\sqrt{h(T)}$ is superfluous, it has not been possible to prove this so far. However, there are other ways in which this problem can be, if not remedied, then at least circumvented, simply by assuming one of the other two mixing conditions.

THEOREM 2.D (Other types of mixing). *Suppose that the requirements for Theorem 2.A are met, with the exception that Condition 3x holds in place of Condition 3ρ.*

(a) *If $x$ is $\beta$, then $d_2(\mathcal{L}(\xi\theta_T^{-1}|_J), \mathcal{L}(\eta\theta_T^{-1}|_J))$ has the same order as that stated in Theorem 2.A, except for the term $\sqrt{Th(T)}\breve{\beta}(m)$, which is replaced by the two terms $\sqrt{T/w(T)}\breve{\alpha}(2^D/w(T))$ and $\sqrt{w(T)T}\breve{\beta}(m)$; hence [since $h(T) \geq 1$ was arbitrary],*

$$d_2(\mathcal{L}(\xi\theta_T^{-1}|_J), \mathcal{L}(\eta\theta_T^{-1}|_J))$$
$$= O\bigg(\frac{1}{T^{1/D_2}}, \log^\uparrow\bigg(\frac{T}{w(T)}\bigg)\frac{m^{D_2}+1}{w(T)},$$
$$\log^\uparrow\bigg(\frac{T}{w(T)}\bigg)\breve{\alpha}\bigg(\frac{2^D(2m+1)^{D_2}}{w(T)}\bigg),$$
$$\sqrt{\frac{T}{w(T)}}\breve{\alpha}\bigg(\frac{2^D}{w(T)}\bigg), \sqrt{w(T)T}\breve{\beta}(m)\bigg) \quad \text{for } T \to \infty.$$

(b) *If $x$ is $\varphi$, then $d_2(\mathcal{L}(\xi\theta_T^{-1}|_J), \mathcal{L}(\eta\theta_T^{-1}|_J))$ has the same order as that stated in Theorem 2.A, but the term $\sqrt{Th(T)}\breve{\beta}(m)$ can be replaced by $\sqrt{T/w(T)}\breve{\beta}(m)$; hence, as above,*

$$d_2(\mathcal{L}(\xi\theta_T^{-1}|_J), \mathcal{L}(\eta\theta_T^{-1}|_J))$$
$$= O\bigg(\frac{1}{T^{1/D_2}}, \log^\uparrow\bigg(\frac{T}{w(T)}\bigg)\frac{m^{D_2}+1}{w(T)},$$
$$\log^\uparrow\bigg(\frac{T}{w(T)}\bigg)\breve{\alpha}\bigg(\frac{2^D(2m+1)^{D_2}}{w(T)}\bigg), \sqrt{\frac{T}{w(T)}}\breve{\beta}(m)\bigg)$$
$$\text{for } T \to \infty.$$

REMARK 2.E. Note that in the above theorem, a certain price must be paid for the elimination of $h(T)$ in the term that comes from the mixing condition: In statement (a) we obtain for our upper bound an order which is, in many cases, worse than the corresponding order we get for an optimal choice of $h(T)$ in Theorem 2.A; only for sufficiently high $D_1$ is the upper bound order from Theorem 2.D(a), in general, better. In statement (b) we require a much stronger kind of mixing condition than in Theorems 2.A and 2.D(a).

On the other hand, we do not have to require a strictly stronger mixing condition in statement (a) and we get a strictly better upper bound in statement (b).

EXAMPLE. A typical choice of parameters for illustrating the above mentioned points is given by $\breve{\alpha}(v) = v$, $\breve{\beta}(u) = \frac{1}{u^{2D_2}}$ and $w(T) = T$, whence we immediately get $O(T^{-1/3})$ and $O(T^{-2/3})$ as upper bound orders for the $d_2$-distance under the $\beta$-mixing and $\varphi$-mixing conditions, respectively; solving a little optimization problem yields the order $O(T^{-3/(D_1+6)})$ under the $\rho$-mixing condition, which for $D_1 < 3$ is better and for $D_1 > 3$ is worse than the order under $\beta$-mixing.

2.2.2. *Proofs.* The following simple lemma will be useful.

LEMMA 2.F. *For all* $\mathbf{k}, \mathbf{l}$, *we have*
$$\nu(C_{\mathbf{kl}}) - 2^{D_2-2} \frac{1}{w(T)h(T)} \breve{\alpha}\left(2^{D_2} \frac{1}{w(T)h(T)}\right) \leq p_{\mathbf{kl}} \leq \nu(C_{\mathbf{kl}}).$$

PROOF. The second inequality is immediate, the first one is obtained as
$$\nu(C_{\mathbf{kl}}) - p_{\mathbf{kl}} = \mathbb{E}\xi(C_{\mathbf{kl}}) - \mathbb{P}[\xi(C_{\mathbf{kl}}) \geq 1]$$
$$= \sum_{r=2}^{\infty}(r-1)\mathbb{P}[\xi(C_{\mathbf{kl}}) = r]$$
$$\leq \frac{1}{4}\mathbb{E}[(\xi(C_{\mathbf{kl}}))^2 \mathbb{1}_{\{\xi(C_{\mathbf{kl}}) \geq 2\}}]$$
$$\leq 2^{D_2-2} \frac{1}{w(T)h(T)} \breve{\alpha}\left(2^{D_2} \frac{1}{w(T)h(T)}\right)$$
by the orderliness condition with $v(C_{\mathbf{kl}}) \leq 2^{D_2} \frac{1}{w(T)h(T)}$. □

PROOF OF THEOREM 2.A. We use the notation introduced in Section 2.1; in particular, we write
$$\Xi := \sum_{\mathbf{k},\mathbf{l}} I_{\mathbf{kl}} \delta_{\alpha_{\mathbf{kl}}} \quad \text{and} \quad \mathrm{H} := \sum_{\mathbf{k},\mathbf{l}} U_{\mathbf{kl}} \delta_{\alpha_{\mathbf{kl}}}$$
for the discretized point processes, where $U_{\mathbf{kl}}$ are independent $\mathrm{Po}(p_{\mathbf{kl}})$-variables for $0 \leq \mathbf{k} \leq 2n_1 + 1$, $0 \leq \mathbf{l} \leq 2n_2 + 1$.

The overall $d_2$-distance can now be split up accordingly:
$$\begin{aligned}(2.1) \quad & d_2(\mathcal{L}(\xi\theta_T^{-1}|_J), \mathcal{L}(\eta\theta_T^{-1}|_J)) \\ & \leq d_2(\mathcal{L}(\xi\theta_T^{-1}|_J), \mathcal{L}(\Xi\theta_T^{-1})) \\ & \quad + d_2(\mathcal{L}(\Xi\theta_T^{-1}), \mathcal{L}(\mathrm{H}\theta_T^{-1})) + d_2(\mathcal{L}(\mathrm{H}\theta_T^{-1}), \mathcal{L}(\eta\theta_T^{-1}|_J)).\end{aligned}$$

We first take a look at the discretization errors. For the $\xi$-discretization we can obtain, via the Kantorovich–Rubinstein equation (1.2),
$$\begin{aligned}(2.2) \quad & d_2(\mathcal{L}(\xi\theta_T^{-1}|_J), \mathcal{L}(\Xi\theta_T^{-1})) \\ & \leq \mathbb{E}d_1(\xi\theta_T^{-1}|_J, \Xi\theta_T^{-1}) \\ & = \mathbb{E}[d_1(\xi\theta_T^{-1}|_J, \Xi\theta_T^{-1})\mathbb{1}_{\{\xi\theta_T^{-1}(J)=W^{(T)}\}}] + 1 \cdot \mathbb{P}[\xi\theta_T^{-1}(J) \neq W^{(T)}].\end{aligned}$$

The second summand can easily be estimated as follows:

$$
\begin{aligned}
\mathbb{P}[\xi\theta_T^{-1}(J) \neq W^{(T)}] &= \mathbb{P}\left[\bigcup_{\mathbf{k},\mathbf{l}}\{\xi(C_{\mathbf{kl}}) \geq 2\}\right] \\
&\leq \sum_{\mathbf{k},\mathbf{l}} \mathbb{P}[\xi(C_{\mathbf{kl}}) \geq 2] \\
&\leq \frac{1}{4}\sum_{\mathbf{k},\mathbf{l}} \mathbb{E}[(\xi(C_{\mathbf{kl}}))^2 \mathbb{1}_{\{\xi(C_{\mathbf{kl}}) \geq 2\}}] \\
&\leq 2^{2D+D_2-2}\frac{T}{w(T)}\breve{\alpha}\left(2^{D_2}\frac{1}{w(T)h(T)}\right)
\end{aligned}
\tag{2.3}
$$

by the orderliness condition with $v(C_{\mathbf{kl}}) \leq 2^{D_2}\frac{1}{w(T)h(T)}$.

In order to estimate the first summand in (2.2), we use the representation of the $d_1$-distance given by (1.1). Let $X_1, \ldots, X_{\xi\theta_T^{-1}(J)}$ be the points of $\xi\theta_T^{-1}|_J$ and $Y_1, \ldots, Y_{W^{(T)}}$ the points of $\Xi\theta_T^{-1}$ and suppose w.l.o.g. that they are numbered in an optimal way on $\{\xi\theta_T^{-1}(J) = W^{(T)}\}$, that is, in such a way that $Y_i$ is the centre $\rho_{\mathbf{kl}}$ of the cuboid $R_{\mathbf{kl}}$ which contains $X_i$. Thus, by (1.1), and since in the transition from $\xi$ to $\Xi$ we do not move the points any farther than half a body diagonal of a cuboid $R_{\mathbf{kl}}$,

$$
\begin{aligned}
&d_1(\xi\theta_T^{-1}|_J, \Xi\theta_T^{-1})\mathbb{1}_{\{\xi\theta_T^{-1}(J)=W^{(T)}\}} \\
&= \left(\frac{1}{W^{(T)}}\sum_{i=1}^{W^{(T)}} d_0(X_i, Y_i)\right)\mathbb{1}_{\{\xi\theta_T^{-1}(J)=W^{(T)}\geq 1\}} \\
&\leq \frac{1}{2}\sqrt{D_1\left(\frac{1}{h(T)^{1/D_1}}\right)^2 + D_2\left(\frac{1}{T^{1/D_2}}\right)^2}\mathbb{1}_{\{\xi\theta_T^{-1}(J)=W^{(T)}\geq 1\}} \\
&\leq \frac{1}{2}\left(\frac{\sqrt{D_1}}{h(T)^{1/D_1}} + \frac{\sqrt{D_2}}{T^{1/D_2}}\right),
\end{aligned}
\tag{2.4}
$$

whence we get for the total $\xi$-discretization error

$$
\begin{aligned}
&d_2(\mathcal{L}(\xi\theta_T^{-1}|_J), \mathcal{L}(\Xi\theta_T^{-1})) \\
&\leq \frac{1}{2}\left(\frac{\sqrt{D_1}}{h(T)^{1/D_1}} + \frac{\sqrt{D_2}}{T^{1/D_2}}\right) + 2^{2D+D_2-2}\frac{T}{w(T)}\breve{\alpha}\left(2^{D_2}\frac{1}{w(T)h(T)}\right).
\end{aligned}
$$

Next we consider the discretization error for $\eta$. Let $\mathrm{H}' := \sum_{\mathbf{k},\mathbf{l}}\eta(C_{\mathbf{kl}})\delta_{\alpha_{\mathbf{kl}}}$ and $q_{\mathbf{kl}} := \nu(C_{\mathbf{kl}})$. We split up the error as

$$
\begin{aligned}
&d_2(\mathcal{L}(\mathrm{H}\theta_T^{-1}), \mathcal{L}(\eta\theta_T^{-1}|_J)) \\
&\leq d_2(\mathcal{L}(\mathrm{H}\theta_T^{-1}), \mathcal{L}(\mathrm{H}'\theta_T^{-1})) + d_2(\mathcal{L}(\mathrm{H}'\theta_T^{-1}), \mathcal{L}(\eta\theta_T^{-1}|_J)).
\end{aligned}
\tag{2.5}
$$

The first summand gives us a little more trouble. Since for any two point processes $\xi_1$ and $\xi_2$ on a compact set $\mathcal{X}$ the inequality

$$
\mathbb{E}d_1(\xi_1, \xi_2) = \mathbb{E}(d_1(\xi_1, \xi_2)\mathbb{1}_{\{\xi_1 \neq \xi_2\}}) \leq \mathbb{P}[\xi_1 \neq \xi_2]
$$

holds, it can be seen from (1.2) and the analogue of (1.4) for probability distributions on more general spaces [see Barbour, Holst and Janson (1992), Appendix A.1] that

$$d_2(P,Q) \leq d_{\mathrm{TV}}(P,Q)$$

for any distributions $P$, $Q$ of point processes on $\mathcal{X}$. Hence, by another application of the more general version of (1.4) in the second inequality,

$$
\begin{aligned}
d_2(\mathcal{L}(\mathrm{H}\theta_T^{-1}), \mathcal{L}(\mathrm{H}'\theta_T^{-1})) &\leq d_{\mathrm{TV}}(\mathcal{L}(\mathrm{H}\theta_T^{-1}), \mathcal{L}(\mathrm{H}'\theta_T^{-1})) \\
&\leq \min_{\substack{U^{(1)}_{\mathbf{kl}} \sim \mathrm{Po}(p_{\mathbf{kl}}), \perp\!\!\!\perp \\ U^{(2)}_{\mathbf{kl}} \sim \mathrm{Po}(q_{\mathbf{kl}}), \perp\!\!\!\perp}} \sum_{\mathbf{k},\mathbf{l}} \mathbb{P}[U^{(1)}_{\mathbf{kl}} \neq U^{(2)}_{\mathbf{kl}}] \\
&= \sum_{\mathbf{k},\mathbf{l}} d_{\mathrm{TV}}(\mathrm{Po}(p_{\mathbf{kl}}), \mathrm{Po}(q_{\mathbf{kl}})) \\
&\leq \sum_{\mathbf{k},\mathbf{l}} (q_{\mathbf{kl}} - p_{\mathbf{kl}}) \\
&\leq 2^{2D+D_2-2} \frac{T}{w(T)} \breve{\alpha}\left(2^{D_2} \frac{1}{w(T)h(T)}\right),
\end{aligned}
$$
(2.6)

where the last two inequalities follow from Proposition A.C and Lemma 2.F, respectively. For the second summand in (2.5), we obtain

$$
\begin{aligned}
d_2(\mathcal{L}(\mathrm{H}'\theta_T^{-1}), \mathcal{L}(\eta\theta_T^{-1}|_J)) &\leq \mathbb{E}d_1(\mathrm{H}'\theta_T^{-1}, \eta\theta_T^{-1}|_J) \\
&= \mathbb{E}[d_1(\mathrm{H}'\theta_T^{-1}, \eta\theta_T^{-1}|_J)\mathbb{1}_{\{\mathrm{H}'\theta_T^{-1}(J) = \eta\theta_T^{-1}(J)\}}] \\
&\leq \frac{1}{2}\left(\frac{\sqrt{D_1}}{h(T)^{1/D_1}} + \frac{\sqrt{D_2}}{T^{1/D_2}}\right)
\end{aligned}
$$
(2.7)

by the same argument that was used in (2.4). So, an estimate for the total $\eta$-discretization error is given by

$$
\begin{aligned}
&d_2(\mathcal{L}(\mathrm{H}\theta_T^{-1}), \mathcal{L}(\eta\theta_T^{-1}|_J)) \\
&\leq \frac{1}{2}\left(\frac{\sqrt{D_1}}{h(T)^{1/D_1}} + \frac{\sqrt{D_2}}{T^{1/D_2}}\right) + 2^{2D+D_2-2}\frac{T}{w(T)}\breve{\alpha}\left(2^{D_2}\frac{1}{w(T)h(T)}\right).
\end{aligned}
$$

Last, we look at the remaining term $d_2(\mathcal{L}(\Xi\theta_T^{-1}), \mathcal{L}(\mathrm{H}\theta_T^{-1}))$, which is perfect for the application of a Stein estimate. In the notation of the Appendix we write

$$\Gamma = \{0, 1, \ldots, 2n_1 + 1\}^{D_1} \times \{0, 1, \ldots, 2n_2 + 1\}^{D_2}$$

[accordingly, we write elements of $\Gamma$ as $(\mathbf{i}, \mathbf{j})$, meaning $\mathbf{i} \in \{0, 1, \ldots, 2n_1 + 1\}^{D_1}$, $\mathbf{j} \in \{0, 1, \ldots, 2n_2 + 1\}^{D_2}$], and for the sets of strongly and weakly dependent indicators, respectively,

$$
\begin{aligned}
\Gamma^s_{\mathbf{kl}} &= \{(\mathbf{i},\mathbf{j}) \in \Gamma_{\mathbf{kl}}; |\mathbf{j} - \mathbf{l}| \leq m\}, \\
\Gamma^w_{\mathbf{kl}} &= \{(\mathbf{i},\mathbf{j}) \in \Gamma_{\mathbf{kl}}; |\mathbf{j} - \mathbf{l}| \geq m + 1\},
\end{aligned}
$$

for every $\mathbf{k}$, $\mathbf{l}$, where $|\mathbf{j} - \mathbf{l}| := \max_{1 \leq s \leq D_2} |j_s - l_s|$ and $m := m(T) \in \mathbb{Z}_+$ for every $T \geq 1$ is chosen arbitrarily. We can assume w.l.o.g. that $m \leq 2n_2 + 1$ [note that for $m > 2n_2 + 1$ we have $e_{\mathbf{kl}} = 0$, so that (2.9) below is still true]. As in the Appendix, we set

$$Z_{\mathbf{kl}} := \sum_{(\mathbf{i},\mathbf{j}) \in \Gamma^s_{\mathbf{kl}}} I_{\mathbf{ij}}, \qquad Y_{\mathbf{kl}} := \sum_{(\mathbf{i},\mathbf{j}) \in \Gamma^w_{\mathbf{kl}}} I_{\mathbf{ij}}.$$

From the local Stein Theorem A.D for point processes we know that

$$
\begin{aligned}
&d_2(\mathcal{L}(\Xi \theta_T^{-1}), \mathcal{L}(\mathrm{H}\theta_T^{-1})) \\
&\quad \leq \left\{ 1 \wedge \frac{2}{\lambda}\left(1 + 2\log^+\left(\frac{\lambda}{2}\right)\right) \right\} \sum_{\mathbf{k},\mathbf{l}} (p_{\mathbf{kl}}^2 + p_{\mathbf{kl}}\mathbb{E}Z_{\mathbf{kl}} + \mathbb{E}(I_{\mathbf{kl}}Z_{\mathbf{kl}})) \\
&\quad + \left(1 \wedge 1.65 \frac{1}{\sqrt{\lambda}}\right) \sum_{\mathbf{k},\mathbf{l}} e_{\mathbf{kl}},
\end{aligned}
$$
(2.8)

with

$$e_{\mathbf{kl}} = 2 \max_{B \in \sigma(I_{\mathbf{ij}};(\mathbf{i},\mathbf{j}) \in \Gamma^w_{\mathbf{kl}})} |\operatorname{cov}(I_{\mathbf{kl}}, \mathbb{1}_B)|.$$

Starting from the right-hand side, most further estimates are very easy. First, we have

$$p_{\mathbf{kl}} \leq \nu(C_{\mathbf{kl}}) \leq \kappa_T \frac{1}{w(T)h(T)}$$

and

$$\mathbb{E}Z_{\mathbf{kl}} = \sum_{\mathbf{i}=0}^{2n_1+1} \sum_{\substack{\mathbf{j}=(\mathbf{l}-m)\vee 0 \\ (\mathbf{i},\mathbf{j}) \neq (\mathbf{k},\mathbf{l})}}^{(\mathbf{l}+m)\wedge(2n_2+1)} p_{\mathbf{ij}} \leq \kappa_T[(2n_1+2)^{D_1}(2m+1)^{D_2} - 1]\frac{1}{w(T)h(T)};$$

furthermore, by the mixing condition,

$$
\begin{aligned}
e_{\mathbf{kl}} &= 2\sqrt{p_{\mathbf{kl}}(1 - p_{\mathbf{kl}})} \max_{B \in \sigma(I_{\mathbf{ij}};(\mathbf{i},\mathbf{j}) \in \Gamma^w_{\mathbf{kl}})} \sqrt{\mathbb{P}[B](1 - \mathbb{P}[B])} |\operatorname{corr}(I_{\mathbf{kl}}, \mathbb{1}_B)| \\
&\leq 2\sqrt{p_{\mathbf{kl}}} \frac{1}{2}\breve{\beta}(m) \leq \sqrt{\kappa_T} \sqrt{\frac{1}{w(T)h(T)}} \breve{\beta}(m);
\end{aligned}
$$
(2.9)

and, by Lemma 2.F,

$$
\begin{aligned}
\lambda &= \sum_{\mathbf{k},\mathbf{l}} p_{\mathbf{kl}} \geq \sum_{\mathbf{k},\mathbf{l}} \left( \nu(C_{\mathbf{kl}}) - 2^{D_2-2} \frac{1}{w(T)h(T)} \breve{\alpha}\left(2^{D_2} \frac{1}{w(T)h(T)}\right) \right) \vee 0 \\
&= \left( \nu(J_T) - (2n_1+2)^{D_1}(2n_2+2)^{D_2} \frac{2^{D_2-2}}{w(T)h(T)} \breve{\alpha}\left(2^{D_2} \frac{1}{w(T)h(T)}\right) \right) \vee 0 \\
&\geq 2^D \frac{T}{w(T)} \left( \iota_T - 2^{D+D_2-2} \breve{\alpha}\left(2^{D_2} \frac{1}{w(T)h(T)}\right) \right) \vee 0,
\end{aligned}
$$

whence we get a "magic factor" estimate of

$$\frac{1}{\lambda} \leq (1 + \varepsilon(T)) \frac{1}{2^D \iota_T} \frac{w(T)}{T},$$

with

$$\varepsilon(T) := \begin{cases} \left(1 - 2^{D+D_2-2} \frac{1}{\iota_T} \check{\alpha}\left(2^{D_2} \frac{1}{w(T)h(T)}\right)\right)^{-1} - 1, & \text{if } (1 - \cdots) > 0, \\ \infty, & \text{otherwise,} \end{cases}$$

an expression of order $O(\check{\alpha}(2^{D_2} \frac{1}{w(T)h(T)}))$ for $T \to \infty$, provided that $\iota > 0$.

For the remaining term, $\mathbb{E}(I_{\mathbf{kl}} Z_{\mathbf{kl}})$, a little trick is required. We subdivide the set $\Gamma = \{0, 1, \ldots, 2n_1 + 1\}^{D_1} \times \{0, 1, \ldots, 2n_2 + 1\}^{D_2}$ along the last $D_2$ dimensions in $D_2$-cube sections of extension $2m+1$ in every dimension (except for possible left over cuboids), and look at the individual sections separately. For $\mathbf{s} = (s_1, s_2, \ldots, s_{D_2}) \in \{1, 2, \ldots, \lceil \frac{2n_2+2}{2m+1} \rceil\}^{D_2}$, set for the $\mathbf{s}$th section, that is, the section containing the $s_j$th collection of $2m + 1$ numbers in the $j$th coordinate,

$$\mathbf{c}^{(1)}(\mathbf{s}) := \mathbf{c}^{(1)}(\mathbf{s}, m) := (c_1^{(1)}(\mathbf{s}), \ldots, c_{D_2}^{(1)}(\mathbf{s}))$$
$$:= ((s_1 - 1)(2m + 1), \ldots, (s_{D_2} - 1)(2m + 1)),$$

which is the "lower left" corner index (the multi-index that is in each coordinate minimal among all indices belonging to the $\mathbf{s}$th section), and

$$\mathbf{c}^{(2)}(\mathbf{s}) := \mathbf{c}^{(2)}(\mathbf{s}, m) := (c_1^{(2)}(\mathbf{s}), \ldots, c_{D_2}^{(2)}(\mathbf{s}))$$
$$:= ([s_1(2m + 1) - 1] \wedge (2n_2 + 1), \ldots, [s_{D_2}(2m + 1) - 1] \wedge (2n_2 + 1)),$$

which is the "upper right" corner index (the multi-index that is in each coordinate maximal among all indices belonging to the $\mathbf{s}$th section). Furthermore, we set

$$D_{\mathbf{s}} := D_{\mathbf{s}}^{(m)} := \bigcup_{\mathbf{i}=0}^{2n_1+1} \bigcup_{\mathbf{j} = [\mathbf{c}^{(1)}(\mathbf{s}) - m] \vee 0}^{[\mathbf{c}^{(2)}(\mathbf{s}) + m] \wedge (2n_2 + 1)} C_{\mathbf{ij}},$$

the subset of $J_T$ that naturally belongs to the $m$-neighborhood cube of the $\mathbf{s}$th section. Using our usual multi-index notation and index range convention for sums, we now obtain for the remaining term

$$\sum_{\mathbf{k},\mathbf{l}} \mathbb{E}(I_{\mathbf{kl}} Z_{\mathbf{kl}})$$

$$= \mathbb{E}\left(\sum_{\mathbf{k}=0}^{2n_1+1} \sum_{\mathbf{l}=0}^{2n_2+1} \sum_{\mathbf{i}=0}^{2n_1+1} \sum_{\substack{\mathbf{j}=(\mathbf{l}-m)\vee 0 \\ (\mathbf{i},\mathbf{j}) \neq (\mathbf{k},\mathbf{l})}}^{(\mathbf{l}+m) \wedge (2n_2+1)} I_{\mathbf{kl}} I_{\mathbf{ij}}\right)$$

BOUNDS FOR POINT PROCESS APPROXIMATIONS 17

$$\leq \mathbb{E}\left\{\sum_{\mathbf{s}=1}^{\lceil(2n_2+2)/(2m+1)\rceil} \left(\sum_{\mathbf{k}=0}^{2n_1+1}\sum_{\mathbf{l}=\mathbf{c}^{(1)}(\mathbf{s})}^{\mathbf{c}^{(2)}(\mathbf{s})}\sum_{\mathbf{i}=0}^{2n_1+1}\sum_{\substack{\mathbf{j}=[\mathbf{c}^{(1)}(\mathbf{s})-m]\vee 0 \\ (\mathbf{i},\mathbf{j})\neq(\mathbf{k},\mathbf{l})}}^{[\mathbf{c}^{(2)}(\mathbf{s})+m]\wedge(2n_2+1)} I_{\mathbf{k}\mathbf{l}}I_{\mathbf{i}\mathbf{j}}\right)\mathbb{1}_{\{\xi(D_\mathbf{s}^{(m)})\geq 2\}}\right\}$$

$$\leq \mathbb{E}\left\{\sum_{\mathbf{s}=1}^{\lceil(2n_2+2)/(2m+1)\rceil}\left(\sum_{\mathbf{i}=0}^{2n_1+1}\sum_{\mathbf{j}=[\mathbf{c}^{(1)}(\mathbf{s})-m]\vee 0}^{[\mathbf{c}^{(2)}(\mathbf{s})+m]\wedge(2n_2+1)} I_{\mathbf{i}\mathbf{j}}\right)^2\mathbb{1}_{\{\xi(D_\mathbf{s}^{(m)})\geq 2\}}\right\}$$

$$\leq \sum_{\mathbf{s}=1}^{\lceil(2n_2+2)/(2m+1)\rceil} \mathbb{E}[(\xi(D_\mathbf{s}^{(m)}))^2 \mathbb{1}_{\{\xi(D_\mathbf{s}^{(m)})\geq 2\}}]$$

$$\leq 2^{D+D_2}(T^{1/D_2}+m+1)^{D_2}\frac{1}{w(T)}\breve{\alpha}\left(2^D(2m+1)^{D_2}\frac{1}{w(T)}\right)$$

by the orderliness condition with $v(D_\mathbf{s}^{(m)}) \leq 2^D(2m+1)^{D_2}\frac{1}{w(T)}$.

All that is left to do now is to combine the various estimates for the right-hand side terms of the Stein inequality (2.8). Then, adding the discretization errors and setting

$$L(T) := 1 \wedge \left[2(1+\varepsilon(T))\frac{w(T)}{2^D \iota_T T}\right]\left(1+2\log^+\left(2^{D-1}\kappa_T \frac{T}{w(T)}\right)\right)$$

yields for the overall $d_2$-distance

(2.10)
$$\begin{aligned}d_2(\mathcal{L}(\xi\theta_T^{-1}|_J),\mathcal{L}(\eta\theta_T^{-1}|_J))\\ \leq \frac{\sqrt{D_1}}{h(T)^{1/D_1}} + \frac{\sqrt{D_2}}{T^{1/D_2}} + L(T)2^{2D+2D_1}\kappa_T^2\frac{T(2m+1)^{D_2}}{(w(T))^2}\\ + 2^{2D+D_2-1}\frac{T}{w(T)}\breve{\alpha}\left(\frac{2^{D_2}}{w(T)h(T)}\right)\\ + L(T)2^{D+D_2}\frac{(T^{1/D_2}+m+1)^{D_2}}{w(T)}\breve{\alpha}\left(2^D\frac{(2m+1)^{D_2}}{w(T)}\right)\\ + \left(1\wedge 1.65\sqrt{1+\varepsilon(T)}\sqrt{\frac{w(T)}{2^D \iota_T T}}\right)2^{2D}\sqrt{\kappa_T}\sqrt{\frac{h(T)}{w(T)}}T\breve{\beta}(m).\end{aligned}$$

For $\iota > 0$ and preferably $T$ large enough, we get the rougher, but less nasty looking upper bound

$$
\begin{aligned}
d_2(&\mathcal{L}(\xi\theta_T^{-1}|_J), \mathcal{L}(\eta\theta_T^{-1}|_J)) \\
&\leq \frac{\sqrt{D_1}}{h(T)^{1/D_1}} + \frac{\sqrt{D_2}}{T^{1/D_2}} \\
&\quad + 2^{D+2D_1+2} \frac{\kappa^2}{\iota}(1+\varepsilon(T))\log^\uparrow\!\left(2^{D-1}\kappa\frac{T}{w(T)}\right)\frac{(2m+1)^{D_2}}{w(T)} \\
&\quad + 2^{2D+D_2-1} \frac{T}{w(T)}\check{\alpha}\!\left(\frac{2^{D_2}}{w(T)h(T)}\right) \\
&\quad + 2^{D_2+2}5^{D_2}\frac{1}{\iota}(1+\varepsilon(T))\log^\uparrow\!\left(2^{D-1}\kappa\frac{T}{w(T)}\right)\check{\alpha}\!\left(2^D\frac{(2m+1)^{D_2}}{w(T)}\right) \\
&\quad + 2^{\frac{3}{2}D+1}\sqrt{\frac{\kappa}{\iota}}\sqrt{1+\varepsilon(T)}\sqrt{Th(T)}\check{\beta}(m),
\end{aligned}
\tag{2.11}
$$

which is of the required order. $\square$

PROOF OF COROLLARY 2.B. For $T \geq 1$, we have to find $h(T) \geq 1$ and $m := m(T) \in \mathbb{Z}_+$, such that all six terms on the right-hand side of the equality in Theorem 2.A go to zero as $T \to \infty$. We set $h(T) = T^q$ and $m := [T^x]$, with $q > 0$ and $0 \leq x < \frac{\delta}{D_2}$. Thus,

$$\frac{1}{h(T)^{1/D_1}} \to 0, \qquad \frac{1}{T^{1/D_2}} \to 0,$$

$$\log^\uparrow\!\left(\frac{T}{w(T)}\right)\frac{m^{D_2}+1}{w(T)} \to 0 \quad \text{and}$$

$$\log^\uparrow\!\left(\frac{T}{w(T)}\right)\check{\alpha}\!\left(\frac{2^D(2m+1)^{D_2}}{w(T)}\right) \to 0;$$

so the only two terms we have to worry about are

$$\frac{T}{w(T)}\check{\alpha}\!\left(\frac{2^{D_2}}{w(T)h(T)}\right) = O(T^{1-\delta-\delta r-qr})$$

and

$$\sqrt{Th(T)}\check{\beta}(m) = O(T^{1/2(1+q-(1+s)D_2x)}),$$

which both converge to zero if there exist $q > 0$ and $0 \leq x < \frac{\delta}{D_2}$ such that

$$q > \frac{1-\delta-\delta r}{r} \quad \text{and} \quad q < (1+s)D_2 x - 1.$$

This last is true provided that

$$(1+s)\delta - 1 > \max\!\left(\frac{1-\delta-\delta r}{r}, 0\right),$$

whence we obtain the statement. □

PROOF OF THEOREM 2.D. Since the mixing condition is used only once in the proof of Theorem 2.A, namely, in (2.9) for obtaining the upper bound of the $e_{\mathbf{kl}}$ from the Stein estimate, we can simply transfer the proof and re-calculate this upper bound under our new mixing conditions.

(a) Let $\mathbf{l} \in \{0, 1, \ldots, 2n_2 + 1\}^{D_2}$ be fixed, set $C_{\cdot \mathbf{l}} := \bigcup_{\mathbf{k}=0}^{2n_1+1} C_{\mathbf{kl}}$, and define

$$\tilde{X}^{(1)}_{\mathrm{int}} := (I_{\mathbf{il}}; \mathbf{i} \in \{0, 1, \ldots, 2n_1 + 1\}^{D_1}), \qquad \tilde{\mathcal{F}}^{(1)}_{\mathrm{int}} := \sigma(\tilde{X}^{(1)}_{\mathrm{int}}),$$
$$\tilde{X}^{(1)}_{\mathrm{ext}} := (I_{\mathbf{ij}}; (\mathbf{i}, \mathbf{j}) \in \Gamma^w_{\mathbf{kl}}) \qquad \text{regardless of } \mathbf{k}, \qquad \tilde{\mathcal{F}}^{(1)}_{\mathrm{ext}} := \sigma(\tilde{X}^{(1)}_{\mathrm{ext}}).$$

Note that $\tilde{\mathcal{F}}^{(1)}_{\mathrm{int}} \subset \mathcal{F}^{(1)}_{\mathrm{int}} := \sigma(\xi|_{C_{\cdot \mathbf{l}}})$ and $\tilde{\mathcal{F}}^{(1)}_{\mathrm{ext}} \subset \mathcal{F}^{(1)}_{\mathrm{ext}} := \sigma(\xi|_{\bigcup_{(\mathbf{i},\mathbf{j}) \in \Gamma^w_{\mathbf{kl}}} C_{\mathbf{ij}}})$, regardless of $\mathbf{k}$. It is seen for every $\mathbf{k} \in \{0, 1, \ldots, 2n_1 + 1\}^{D_1}$ that

$$e_{\mathbf{kl}} = 2 \max_{B \in \tilde{\mathcal{F}}^{(1)}_{\mathrm{ext}}} |\mathrm{cov}(I_{\mathbf{kl}}, \mathbb{1}_B)|$$
$$= 2 \max_{B \in \tilde{\mathcal{F}}^{(1)}_{\mathrm{ext}}} |\mathbb{P}[B \cap \{I_{\mathbf{kl}} = 1\}] - \mathbb{P}[B]\mathbb{P}[I_{\mathbf{kl}} = 1]|$$
$$\leq 2 \max_{B \in \tilde{\mathcal{F}}^{(1)}_{\mathrm{ext}}} |\mathbb{P}[B \cap \{\tilde{X}^{(1)}_{\mathrm{int}} = x_{\mathbf{k}}\}] - \mathbb{P}[B]\mathbb{P}[\tilde{X}^{(1)}_{\mathrm{int}} = x_{\mathbf{k}}]|$$
$$+ 2 \max_{B \in \tilde{\mathcal{F}}^{(1)}_{\mathrm{ext}}} \left| \mathbb{P}\left[B \cap \{I_{\mathbf{kl}} = 1\} \cap \left\{\sum_{\mathbf{i}} I_{\mathbf{il}} \geq 2\right\}\right] \right.$$
$$\left. - \mathbb{P}[B]\mathbb{P}\left[\{I_{\mathbf{kl}} = 1\} \cap \left\{\sum_{\mathbf{i}} I_{\mathbf{il}} \geq 2\right\}\right] \right|,$$

where $x_{\mathbf{k}}$ is the element of $\{0, 1\}^{\{0, 1, \ldots, 2n_1+1\}^{D_1}}$, which has a 1 in the $\mathbf{k}$th and a 0 in every other component. We denote the first summand by $A_{\mathbf{kl}}$, the second by $B_{\mathbf{kl}}$ and look at the sums over $\mathbf{k}$ separately. For the $A_{\mathbf{kl}}$-sum we obtain

$$\sum_{\mathbf{k}=0}^{2n_1+1} A_{\mathbf{kl}} = 2 \sum_{\mathbf{k}=0}^{2n_1+1} \max_{B \in \tilde{\mathcal{F}}^{(1)}_{\mathrm{ext}}} |\mathbb{P}[B|\tilde{X}^{(1)}_{\mathrm{int}} = x_{\mathbf{k}}] - \mathbb{P}[B]|\mathbb{P}[\tilde{X}^{(1)}_{\mathrm{int}} = x_{\mathbf{k}}]$$
$$\leq 2\mathbb{E}\left(\max_{B \in \tilde{\mathcal{F}}^{(1)}_{\mathrm{ext}}} |\mathbb{P}[B|\tilde{X}^{(1)}_{\mathrm{int}}] - \mathbb{P}[B]|\right)$$
$$= 2\beta(\tilde{\mathcal{F}}^{(1)}_{\mathrm{int}}, \tilde{\mathcal{F}}^{(1)}_{\mathrm{ext}}) \leq 2\beta(\mathcal{F}^{(1)}_{\mathrm{int}}, \mathcal{F}^{(1)}_{\mathrm{ext}}) \leq 2\breve{\beta}(m),$$

where the monotony of the $\beta$-mixing coefficient is immediate if it is written in its dual form as a supremum over measurable partitions [see Doukhan (1994),

Section 1.1]. For the $B_{\mathbf{kl}}$-sum, the upper bound is obtained by application of the orderliness condition:

$$\sum_{\mathbf{k}=0}^{2n_1+1} B_{\mathbf{kl}} \leq 4 \sum_{\mathbf{k}=0}^{2n_1+1} \mathbb{E}(I_{\mathbf{kl}}\mathbb{1}_{\{\sum_{\mathbf{i}} I_{\mathbf{il}} \geq 2\}})$$

$$\leq 2\mathbb{E}[(\xi(C_{\cdot \mathbf{l}}))^2 \mathbb{1}_{\{\xi(C_{\cdot \mathbf{l}}) \geq 2\}}]$$

$$\leq 2^{D+1} \frac{1}{w(T)} \breve{\alpha}\left(2^D \frac{1}{w(T)}\right).$$

We thus have for the total $e_{\mathbf{kl}}$-sum over $\mathbf{k}$ the estimate

$$\sum_{\mathbf{k}=0}^{2n_1+1} e_{\mathbf{kl}} \leq 2\breve{\beta}(m) + 2^{D+1} \frac{1}{w(T)} \breve{\alpha}\left(2^D \frac{1}{w(T)}\right).$$

(b) In the case of the $\varphi$-mixing condition, the corresponding estimate is very easy. It follows that

$$e_{\mathbf{kl}} = 2 \max_{B \in \tilde{\mathcal{F}}_{\text{ext}}^{(\mathbf{l})}} |\text{cov}(I_{\mathbf{kl}}, \mathbb{1}_B)|$$

$$= 2\left(\max_{B \in \tilde{\mathcal{F}}_{\text{ext}}^{(\mathbf{l})}} |\mathbb{P}[B|I_{\mathbf{kl}} = 1] - \mathbb{P}[B]|\right)\mathbb{P}[I_{\mathbf{kl}} = 1]$$

$$\leq 2\breve{\beta}(m) \frac{\kappa_T}{h(T)w(T)}. \qquad \square$$

2.3. *The $d_{\text{TV}}$-distance between the numbers of points.* Since for every $A \subset \mathbb{Z}_+$ the function $f_A \colon \mathcal{M}_p \to \mathbb{R}_+$ that is defined by $f_A(\rho) := \mathrm{I}[|\rho| \in A]$ is in $\mathcal{F}_2$, it follows for any two point processes $\xi_1, \xi_2$ on a compact set $\mathcal{X}$, that

$$|\mathbb{P}[\xi_1(\mathcal{X}) \in A] - \mathbb{P}[\xi_2(\mathcal{X}) \in A]| \leq d_2(\mathcal{L}(\xi_1), \mathcal{L}(\xi_2)),$$

hence, also

$$d_{\text{TV}}(\mathcal{L}(\xi_1(\mathcal{X})), \mathcal{L}(\xi_2(\mathcal{X}))) \leq d_2(\mathcal{L}(\xi_1), \mathcal{L}(\xi_2)).$$

Thus, the upper bounds we obtained in the theorems of Section 2.2 are also upper bounds for $d_{\text{TV}}(\mathcal{L}(\xi \theta_T^{-1}(J)), \mathcal{L}(\eta \theta_T^{-1}(J)))$. However, using the same method as above and making only slight modifications in the proofs, one can do a little better. Note that although now we are only concerned about numbers of points and not about their positions, we can still improve (but possibly also impair, depending on the leading term in our estimate) our upper bound by choosing a finer discretization in the $\mathbb{R}^{D_1}$-directions. This is because the advantage we get from the orderliness condition if we have smaller discretization cuboids surmounts the disadvantage of having more of them.

THEOREM 2.G. *Suppose that the prerequisites of Section 1 hold, including the Conditions 1, 2 and 3$\rho$, and let $\iota > 0$.*

*Then we obtain for arbitrary $m := m(T) \in \mathbb{Z}_+$ and $h(T) \geq 1$ for every $T \geq 1$:*

$$d_{\mathrm{TV}}(\mathcal{L}(\xi\theta_T^{-1}(J)), \mathcal{L}(\eta\theta_T^{-1}(J)))$$
$$= O\left(\frac{m^{D_2}+1}{w(T)}, \frac{T}{w(T)}\breve{\alpha}\left(\frac{2^{D_2}}{w(T)h(T)}\right), \breve{\alpha}\left(\frac{2^D(2m+1)^{D_2}}{w(T)}\right), \sqrt{Th(T)}\breve{\beta}(m)\right)$$
$$\text{for } T \to \infty.$$

REMARK 2.H. Of course, all theorems stated in Section 2.2 have their equivalents for the $d_{\mathrm{TV}}$-distance between the distributions of the numbers of points. The corresponding upper bounds can simply be obtained by leaving out the $\log^{\uparrow}$-terms, as well as the terms

$$\frac{1}{h(T)^{1/D_1}} \quad \text{and} \quad \frac{1}{T^{1/D_2}}.$$

Note, however, that the conditions in Corollary 2.B for convergence to zero of the principal upper bound remain unchanged.

PROOF OF THEOREM 2.G. Although our task now seems to be quite different, we can proceed exactly as we did in the proof of Theorem 2.A. First, we split up the distance as

$$\begin{aligned} &d_{\mathrm{TV}}(\mathcal{L}(\xi\theta_T^{-1}(J)), \mathcal{L}(\eta\theta_T^{-1}(J))) \\ &= d_{\mathrm{TV}}(\mathcal{L}(\xi(J_T)), \mathrm{Po}(\nu(J_T))) \\ &\leq d_{\mathrm{TV}}(\mathcal{L}(\xi(J_T)), \mathcal{L}(W)) \\ &\quad + d_{\mathrm{TV}}(\mathcal{L}(W), \mathrm{Po}(\lambda)) + d_{\mathrm{TV}}(\mathrm{Po}(\lambda), \mathrm{Po}(\nu(J_T))). \end{aligned}$$

Here the two discretization errors can be estimated very easily. By the orderliness condition, we obtain

$$\begin{aligned} d_{\mathrm{TV}}(\mathcal{L}(\xi(J_T)), \mathcal{L}(W)) &\leq \mathbb{P}[\xi(J_T) \neq W] \\ &= \mathbb{P}\left[\bigcup_{\mathbf{k},\mathbf{l}}\{\xi(C_{\mathbf{kl}}) \geq 2\}\right] \\ &\leq \frac{1}{4}\sum_{\mathbf{k},\mathbf{l}} \mathbb{E}[(\xi(C_{\mathbf{kl}}))^2 \mathbb{1}_{\{\xi(C_{\mathbf{kl}})\geq 2\}}] \\ &\leq 2^{2D+D_2-2}\frac{T}{w(T)}\breve{\alpha}\left(2^{D_2}\frac{1}{w(T)h(T)}\right) \end{aligned}$$

and by Proposition A.C,

$$d_{\mathrm{TV}}(\mathrm{Po}(\lambda), \mathrm{Po}(\nu(J_T)))$$
$$\leq \min\biggl(1, \frac{1}{\sqrt{\lambda}}, \frac{1}{\sqrt{\nu(J_T)}}\biggr)|\lambda - \nu(J_T)|$$
$$= \biggl(1 \wedge \frac{1}{\sqrt{\nu(J_T)}}\biggr) \sum_{\mathbf{k},\mathbf{l}} (\nu(C_{\mathbf{k}\mathbf{l}}) - p_{\mathbf{k}\mathbf{l}})$$
$$\leq \biggl(1 \wedge \frac{1}{2^{D/2}\sqrt{\iota_T}} \sqrt{\frac{w(T)}{T}}\biggr) 2^{2D+D_2-2} \frac{T}{w(T)} \check{\alpha}\biggl(2^{D_2} \frac{1}{w(T)h(T)}\biggr).$$

As for the remaining term, $d_{\mathrm{TV}}(\mathcal{L}(W), \mathrm{Po}(\lambda))$, we can proceed exactly as we did with $d_2(\mathcal{L}(\Xi\theta_T^{-1}), \mathcal{L}(\mathrm{H}\theta_T^{-1}))$, with the only difference that now we use the classical local Stein–Chen Theorem A.A. Thus,

$$d_{\mathrm{TV}}(\mathcal{L}(W), \mathrm{Po}(\lambda))$$
$$\leq \min\biggl(1, \frac{1}{\lambda}\biggr) \sum_{\mathbf{k},\mathbf{l}} (p_{\mathbf{k}\mathbf{l}}^2 + p_{\mathbf{k}\mathbf{l}} \mathbb{E} Z_{\mathbf{k}\mathbf{l}} + \mathbb{E}(I_{\mathbf{k}\mathbf{l}} Z_{\mathbf{k}\mathbf{l}})) + \min\biggl(1, \frac{1}{\sqrt{\lambda}}\biggr) \sum_{\mathbf{k},\mathbf{l}} e_{\mathbf{k}\mathbf{l}}$$

with

$$e_{\mathbf{k}\mathbf{l}} = 2 \max_{B \in \sigma(I_{\mathbf{i}\mathbf{j}}; (\mathbf{i},\mathbf{j}) \in \Gamma_{\mathbf{k}\mathbf{l}}^w)} |\mathrm{cov}(I_{\mathbf{k}\mathbf{l}}, \mathbb{1}_B)|.$$

All notation has exactly the same meaning as it had in the proof of Theorem 2.A, so except for the logarithmic factor in front of the first sum, and the constant 1.65 in front of the second, we get exactly the same upper bound for $d_{\mathrm{TV}}(\mathcal{L}(W), \mathrm{Po}(\lambda))$ as we did for $d_2(\mathcal{L}(\Xi\theta_T^{-1}), \mathcal{L}(\mathrm{H}\theta_T^{-1}))$.

Assembling of all the different pieces yields the result claimed. □

2.4. *Results for measure preserving transformations* $\tilde{\theta}_T$. When we consider a stretch factor $w(T)^{1/D_1} = o(T^{1/D_1})$, the expected number of points of the transformed process $\xi\theta_T^{-1}$ contained within the fixed cube $J$ goes to infinity as $T \to \infty$ if $\iota > 0$, which for some applications is not desirable (e.g., if we want to approximate $\xi\theta_T^{-1}|_J$ by a Poisson process that does not depend on $T$, see Section 2.5). We therefore formulate another theorem in this section, which deals with the case where we adjust the volume of the cuboid $J$ to the volume of the cuboids $J_T$, and thus produce space for the additional points.

In this regard, let $\tilde{\theta}_T$ and $\tilde{J}_T$, defined as in Section 1, be our substitute for the transformation $\theta_T$ and our enlarged version of the cuboid $J$, respectively. We then obtain the following result, where once more the quantitative form of the upper bound can be found at the end of the proof.

THEOREM 2.I. *Suppose that the prerequisites of Section* 1 *hold, including the Conditions* 1, 2 *and* 3$\rho$, *and let* $\iota > 0$.

*Then we obtain for arbitrary* $m := m(T) \in \mathbb{Z}_+$ *and* $h(T) \geq 1$ *for every* $T \geq 1$:

$$d_2(\mathcal{L}(\xi\tilde{\theta}_T^{-1}|_{\tilde{J}_T}), \mathcal{L}(\eta\tilde{\theta}_T^{-1}|_{\tilde{J}_T}))$$
$$= O\bigg(\bigg(\frac{T}{w(T)}\bigg)^{1/D_1}\frac{1}{h(T)^{1/D_1}}, \frac{1}{T^{1/D_2}},$$
$$\log^{\uparrow}\bigg(\frac{T}{w(T)}\bigg)\frac{m^{D_2}+1}{w(T)}, \frac{T}{w(T)}\check{\alpha}\bigg(\frac{2^{D_2}}{w(T)h(T)}\bigg),$$
$$\log^{\uparrow}\bigg(\frac{T}{w(T)}\bigg)\check{\alpha}\bigg(\frac{2^D(2m+1)^{D_2}}{w(T)}\bigg), \sqrt{Th(T)}\check{\beta}(m)\bigg)$$
$$\text{for } T \to \infty,$$

*which is the same order as in Theorem* 2.A, *apart from the factor* $(T/W(T))^{1/D_1}$.

PROOF. For a large part we can adopt the proof of Theorem 2.A. We use the same notation and the same discretization as we did there, replacing only $\theta_T$ by $\tilde{\theta}_T$ and $J$ by $\tilde{J}_T$. First note that there is no change at all for the estimate of the Stein term, now written as $d_2(\mathcal{L}(\Xi\tilde{\theta}_T^{-1}), \mathcal{L}(H\tilde{\theta}_T^{-1}))$, because in the Stein estimate only objects in the pre-image of $\tilde{\theta}_T$ have to be considered (the Stein estimate does not take into account the distances between the points!).

But the changes for the estimates of the approximation errors are not exactly huge either: As we have seen in the proof of Theorem 2.A, these errors can be split up into two additive parts, one stemming from the fact that the original and the discretized point process need not have the same numbers of points in every discretization cuboid [see (2.3), resp. (2.6), in the proof of Theorem 2.A] and one stemming from the fact that even when we have the same numbers of points in every discretization cuboid, their positions are, in general, a bit shifted [see (2.4), resp. (2.7)]. From those two parts only the second is affected by the transition from $\theta_T$ to $\tilde{\theta}_T$ and from $J$ to $\tilde{J}_T$ (inasmuch as the discretization cuboids in the image space get a little bigger), because for the first, we have to deal once more only with objects in the pre-image of $\tilde{\theta}_T$. A short calculation taking into account the above considerations [reproducing inequalities (2.4) and, accordingly, (2.7)] provides as upper bounds for each of the discretization errors $d_2(\mathcal{L}(\xi\tilde{\theta}_T^{-1}|_{\tilde{J}_T}), \mathcal{L}(\Xi\tilde{\theta}_T^{-1}))$ and $d_2(\mathcal{L}(H\tilde{\theta}_T^{-1}), \mathcal{L}(\eta\tilde{\theta}_T^{-1}|_{\tilde{J}_T}))$,

$$\frac{1}{2}\bigg(\bigg(\frac{T}{w(T)}\bigg)^{1/D_1}\frac{\sqrt{D_1}}{h(T)^{1/D_1}} + \frac{\sqrt{D_2}}{T^{1/D_2}}\bigg) + 2^{2D+D_2-2}\frac{T}{w(T)}\check{\alpha}\bigg(2^{D_2}\frac{1}{w(T)h(T)}\bigg).$$

Thus, we obtain as possible upper bounds for the overall $d_2$-distance those of (2.10) and (2.11) with $\frac{\sqrt{D_1}}{h(T)^{1/D_1}}$ replaced by $(\frac{T}{w(T)})^{1/D_1} \frac{\sqrt{D_1}}{h(T)^{1/D_1}}$, which yields the required qualitative estimate. □

Again we can formulate versions of the other results of Section 2.2 with only slight (and very obvious) changes; in particular, we get the following:

COROLLARY 2.J (Convergence to zero in Theorem 2.I). *Suppose that the prerequisites of Theorem 2.I hold. Furthermore, suppose that $w(T) \geq kT^\delta$ for $k > 0$, $\delta \in (0,1]$ and that*

$$\breve{\alpha}(v) = O(v^r) \qquad \text{for } v \to 0 \text{ with } r > 0,$$

$$\breve{\beta}(u) = O\left(\frac{1}{u^{(1+s)D_2/2}}\right) \qquad \text{for } u \to \infty \text{ with } 1+s > \max\left(\frac{1-\delta}{\delta}\frac{1+r}{r}, \frac{2-\delta}{\delta}\right).$$

*Then*

$$d_2(\mathcal{L}(\xi\tilde{\theta}_T^{-1}|_{\tilde{J}_T}), \mathcal{L}(\eta\tilde{\theta}_T^{-1}|_{\tilde{J}_T})) \to 0 \qquad \text{for } T \to \infty.$$

Note that under the $\beta$-mixing or the $\varphi$-mixing condition, no changes in the respective upper bound order obtained in Theorem 2.D are necessary.

2.5. *Results for a fixed limiting process.* So far we have only examined approximations of the transformed process $\xi\theta_T^{-1}$ (resp. $\xi\tilde{\theta}_T^{-1}$) by a Poisson process which has the expectation measure $\nu\theta_T^{-1}$. Of course, this implies that the expectation measure may (and, unless it is a constant multiple of the Lebesgue measure, does) change as $T$ tends to infinity: The approximating Poisson process, in general, will not be stable. One might therefore ask under what circumstances it is possible to approximate the transformed $\xi$-process by a fixed Poisson process, whose distribution does not depend on $T$, and what loss in terms of the $d_2$-distance one has to face.

First of all, the correct $T$-independent intensity measure for our new Poisson process has to be found. Clearly, for $\iota > 0$, using the transformation $\theta_T$ with a stretch factor $w(T) = o(T)$ is unnatural, because in that case the expected number of points of $\xi\theta_T^{-1}$ contained in $J$ goes to infinity, whereas, of course, for any fixed Poisson process, the expectation of the number of points in $J$ is always finite. So the natural choice for general $w(T)$ is the measure preserving transformation $\tilde{\theta}_T$, together with the enlarged cuboid $\tilde{J}_T$ from Section 2.4.

For the following heuristics we ignore the fact that $\mu_2$ might be a counting measure. Then, restricted to the cuboid $J_T$ for $T$ relatively large, the measure $\nu$ with density $p$ with respect to $\lambda^D$ should be relatively "close" to the measure $\nu' := p(\mathbf{0})\lambda^D$, provided that $p$ is constant in the $\mathbb{R}^{D_2}$-directions [hence, the notation $p(\mathbf{s}) = p(\mathbf{s}, \mathbf{t})$ for all $\mathbf{s} \in \mathbb{R}^{D_1}$, $\mathbf{t} \in \mathbb{R}^{D_2}$] and that $p$ satisfies a regularity condition in the $\mathbb{R}^{D_1}$-directions at $\mathbf{0}$. Thus, restricted to $\tilde{J}_T$, $\nu \tilde{\theta}_T^{-1}$ should be close to $\nu' \tilde{\theta}_T^{-1}$ [which is again $p(\mathbf{0})\lambda^D$, hence, not dependent on $T$] as well, and, therefore, $\mathrm{Po}(p(\mathbf{0})\lambda^D|_{\tilde{J}_T})$ should be a good choice for approximating $\mathcal{L}(\xi \theta_T^{-1}|_{\tilde{J}_T})$.

The following makes the above considerations rigorous. First, we formulate the additional regularity condition for $p$.

CONDITION 4 (Regularity of $p$). The density $p = d\nu/d\mu$ is constant in the $\mathbb{R}^{D_2}$-directions, so that we can write

$$p(\mathbf{s}, \mathbf{t}) = p(\mathbf{s}) \qquad \text{for all } \mathbf{s} \in \mathbb{R}^{D_1}, \mathbf{t} \in \mathbb{R}^{D_2} (\text{resp. } \mathbf{t} \in \mathbb{Z}^{D_2} + \tfrac{1}{2}\mathbf{1}).$$

Moreover, $p$ satisfies the following regularity condition in the $\mathbb{R}^{D_1}$-directions: There exist $L \geq 0$ and $z > 0$, such that

$$|p(\mathbf{s}) - p(\mathbf{0})| \leq L|\mathbf{s}|^z \qquad \text{for all } \mathbf{s} \in \mathbb{R}^{D_1}$$

(or for $\mathbf{s} \in [-(\frac{1}{w(T)})^{1/D_1}, (\frac{1}{w(T)})^{1/D_1}]^{D_1}$ for the $T$ one wishes to consider).

We are now in the position to formulate the theorem.

THEOREM 2.K. *Suppose that the prerequisites of Section* 1 *hold, including the Conditions* 1, 2, 3$\rho$, *as well as the new Condition* 4 *above. Let* $\iota > 0$, $T \geq 1$ *(remember that we always assume that* $T \in \{n^{D_2}; n \in \mathbb{N}\}$ *if* $\mu_2 = \mathcal{H}_0^{D_2}$*),* $m := m(T) \in \mathbb{Z}_+$, *and* $h(T) \geq 1$. *Then*

$$d_2(\mathcal{L}(\xi \tilde{\theta}_T^{-1}|_{\tilde{J}_T}), \mathrm{Po}(p(\mathbf{0})\lambda^D|_{\tilde{J}_T}))$$

$$\leq \tilde{A}(T) + 2^{(z+D_1+2D_2)/2} \frac{D_1}{z+D_1} L \tau_{D_1} \frac{T}{w(T)^{1+z/D_1}}$$

$$= O\left(\frac{T}{w(T)^{1+z/D_1}}, \left(\frac{T}{w(T)}\right)^{1/D_1} \frac{1}{h(T)^{1/D_1}}, \frac{1}{T^{1/D_2}},\right.$$

$$\log^{\uparrow}\left(\frac{T}{w(T)}\right) \frac{m^{D_2}+1}{w(T)}, \frac{T}{w(T)} \breve{\alpha}\left(\frac{2^{D_2}}{w(T)h(T)}\right),$$

$$\left.\log^{\uparrow}\left(\frac{T}{w(T)}\right) \breve{\alpha}\left(\frac{2^D(2m+1)^{D_2}}{w(T)}\right), \sqrt{Th(T)} \breve{\beta}(m)\right)$$

$$\text{for } T \to \infty,$$

*where* $\tilde{A}(T) := \tilde{A}(T, m, h(T))$ *is the explicit upper bound that we obtained in Theorem* 2.I [*formula* (2.10) *or* (2.11) *with the corresponding modifications*] *and* $\tau_{D_1} = \pi^{D_1/2}/\Gamma(\frac{D_1}{2}+1)$ *is the volume of the $D_1$-dimensional unit ball.*

COROLLARY 2.L. *Under the prerequisites of Corollary 2.J plus Condition 4, with $z > \frac{1-\delta}{\delta} D_1$, we obtain*

$$d_2(\mathcal{L}(\xi \tilde{\theta}_T^{-1}|_{\tilde{J}_T}), \operatorname{Po}(p(\mathbf{0}) \lambda^D|_{\tilde{J}_T})) \to 0 \qquad \text{for } T \to \infty,$$

*hence, if $\delta = 1$ ($z > 0$),*

$$\xi \tilde{\theta}_T^{-1}|_J \xrightarrow{\mathcal{D}} \operatorname{Po}(p(\mathbf{0}) \lambda^D|_J),$$

*by result* (1.3).

PROOF OF THEOREM 2.K. Once again we can largely adopt the proof of Theorem 2.A (or, more precisely, that of Theorem 2.I). This time only the estimate for the discretization error $d_2(\mathcal{L}(\mathrm{H}\tilde{\theta}_T^{-1}), \mathcal{L}(\eta \tilde{\theta}_T^{-1}|_{\tilde{J}_T}))$ has to be replaced by an appropriate estimate for our new error $d_2(\mathcal{L}(\mathrm{H}\tilde{\theta}_T^{-1}), \operatorname{Po}(p(\mathbf{0})\lambda^D|_{\tilde{J}_T}))$. We proceed just as we did in Theorem 2.A.

Let $\eta' \sim \operatorname{Po}(p(\mathbf{0})\lambda^D)$ [consequently, also $\eta' \tilde{\theta}_T^{-1} \sim \operatorname{Po}(p(\mathbf{0})\lambda^D)$], $\mathrm{H}'' := \sum_{\mathbf{k},\mathbf{l}} \eta'(C_{\mathbf{k}\mathbf{l}}) \delta_{\alpha_{\mathbf{k}\mathbf{l}}}$, and split up the error as

$$\begin{aligned} d_2&(\mathcal{L}(\mathrm{H}\tilde{\theta}_T^{-1}), \operatorname{Po}(p(\mathbf{0})\lambda^D|_{\tilde{J}_T})) \\ &= d_2(\mathcal{L}(\mathrm{H}\tilde{\theta}_T^{-1}), \mathcal{L}(\eta' \tilde{\theta}_T^{-1}|_{\tilde{J}_T})) \\ &\leq d_2(\mathcal{L}(\mathrm{H}\tilde{\theta}_T^{-1}), \mathcal{L}(\mathrm{H}'' \tilde{\theta}_T^{-1})) + d_2(\mathcal{L}(\mathrm{H}'' \tilde{\theta}_T^{-1}), \mathcal{L}(\eta' \tilde{\theta}_T^{-1}|_{\tilde{J}_T})). \end{aligned}$$

Inequality (2.7) (or, more precisely, the corresponding modification from the proof of Theorem 2.I) yields for the second summand, as before,

$$(2.12) \quad d_2(\mathcal{L}(\mathrm{H}'' \tilde{\theta}_T^{-1}), \mathcal{L}(\eta' \tilde{\theta}_T^{-1}|_{\tilde{J}_T})) \leq \frac{1}{2}\left( \left(\frac{T}{w(T)}\right)^{1/D_1} \frac{\sqrt{D_1}}{h(T)^{1/D_1}} + \frac{\sqrt{D_2}}{T^{1/D_2}} \right).$$

For the first summand we get, by the same method as in (2.6),

$$(2.13) \quad \begin{aligned} d_2&(\mathcal{L}(\mathrm{H}\tilde{\theta}_T^{-1}), \mathcal{L}(\mathrm{H}'' \tilde{\theta}_T^{-1})) \\ &\leq \sum_{\mathbf{k},\mathbf{l}} d_{\mathrm{TV}}(\operatorname{Po}(p_{\mathbf{k}\mathbf{l}}), \operatorname{Po}(p(\mathbf{0})\lambda^D(C_{\mathbf{k}\mathbf{l}}))) \\ &\leq \sum_{\mathbf{k},\mathbf{l}} (\nu(C_{\mathbf{k}\mathbf{l}}) - p_{\mathbf{k}\mathbf{l}}) + \sum_{\mathbf{k},\mathbf{l}} |\nu(C_{\mathbf{k}\mathbf{l}}) - p(\mathbf{0})\lambda^D(C_{\mathbf{k}\mathbf{l}})|, \end{aligned}$$

where the first sum was already estimated in (2.6). Its upper bound, together with the upper bound from (2.12), forms the bound we arrived at for $d_2(\mathcal{L}(\mathrm{H}\tilde{\theta}_T^{-1}), \mathcal{L}(\eta \tilde{\theta}_T^{-1}|_{\tilde{J}_T}))$. Therefore, all that is left to do is to show that the second sum on the right-hand side of (2.13) can be estimated by the claimed additional term. This, however, is done very easily:

$$\sum_{\mathbf{k},\mathbf{l}} |\nu(C_{\mathbf{k}\mathbf{l}}) - p(\mathbf{0})\lambda^D(C_{\mathbf{k}\mathbf{l}})| = \sum_{\mathbf{k},\mathbf{l}} \left| \int_{C_{\mathbf{k}\mathbf{l}}} (p(\mathbf{s}) - p(\mathbf{0})) \mu(d(\mathbf{s}, \mathbf{t})) \right|$$

$$\leq \int_{J_T} |p(\mathbf{s}) - p(\mathbf{0})| \mu(d(\mathbf{s},\mathbf{t}))$$

$$\leq 2^{D_2} L \cdot T \int_{[-(1/w(T))^{1/D_1}, (1/w(T))^{1/D_1}]^{D_1}} |\mathbf{s}|^z \lambda^{D_1}(d\mathbf{s})$$

$$\leq 2^{D_2} D_1 L \tau_{D_1} \cdot T \int_0^{\sqrt{2}(1/w(T))^{1/D_1}} r^{z+D_1-1} dr$$

$$= 2^{(z+D_1+2D_2)/2} \frac{D_1}{z+D_1} L \tau_{D_1} \frac{T}{w(T)^{1+z/D_1}}. \qquad \square$$

**3. Applications.** The results of Section 2 can be applied in a number of different ways. For example, they yield useful upper bounds for certain theoretical statements about Poisson process approximation, such as classical thinning and superposition theorems (by projection of the point processes involved on the $\mathbb{R}^{D_2}$-directions and the $\mathbb{R}^{D_1}$-directions, resp.). There are also statistical problems where the results of Section 2 can be of help. To obtain an idea of what is possible, we look at two examples in more detail: in Section 3.1 we consider a fairly general density estimation problem, examined by Ellis (1991), and in Section 3.2 we consider a problem of testing for long range dependence.

3.1. *Density estimation.* First of all, we need a new regularity condition for the density $p$.

CONDITION 4' (Regularity of $p$). The density $p = d\nu/d\mu$ is constant in the $\mathbb{R}^{D_2}$-directions, so that we can write

$$p(\mathbf{s},\mathbf{t}) = p(\mathbf{s}) \qquad \text{for all } \mathbf{s} \in \mathbb{R}^{D_1}, \mathbf{t} \in \mathbb{R}^{D_2} (\text{resp. } \mathbf{t} \in \mathbb{Z}^{D_2} + \tfrac{1}{2}\mathbf{1}).$$

Moreover, $p$ satisfies the following regularity condition in the $\mathbb{R}^{D_1}$-directions:

$$p \in C^2(\mathbb{R}^{D_1}).$$

Of course, it is enough if $p|_Z \in C^2(Z)$ for a sufficiently large neighborhood $Z$ of $\mathbf{0} \in \mathbb{R}^{D_1}$.

Suppose that Condition 4' holds (along with the usual conditions from Section 1), and that we want to estimate the density $p$ at the point $\mathbf{0} \in \mathbb{R}^{D_1}$, say.

By way of illustration, it is convenient to think of the $\mathbb{R}^{D_1}$-space as the "data space" (i.e., the space of possible data points) and the $\mathbb{R}^{D_2}$-space as the "ascertainment space" [i.e., the space of points at which data is obtained, typically by continuous observation over time ($\mathbb{R}^{D_2} = \mathbb{R} = $ time axis) or by repetition of experiments ($\mathbb{R}^{D_2}$ with reference measure $\mu_2 = \mathcal{H}_0^{D_2}$)]. An example suggested by Ellis (1986, 1991) is the estimation of the rate at which earthquakes above a certain magnitude occur per unit area and unit time in a certain region. Here we have $D_1 = 2$ and $D_2 = 1$, and the points in $\mathbb{R}^3$ represent the positions and times of the observed earthquakes.

Among various methods for density estimation, we choose kernel estimation with a data-independent window width, that is, the window width in the $\mathbb{R}^{D_1}$-directions does not depend directly on the data, but does depend on the "observation span" (which in the discrete case corresponds to the sample size). For a detailed account of density estimation see Silverman (1986). We adapt the usual notation in connection with density estimation to the notation we used in Section 2. Thus, $2T^{1/D_2}$ is our observation span (in $D_2$ directions), $2/w(T)^{1/D_1}$ is the window width (in $D_1$ directions) and our density estimator at the point $\mathbf{0}$ takes the form

$$\hat{p}_\xi(\mathbf{0}) := \frac{1}{|J_T|} \int_{J_T} 2^{D_1} K(w(T)^{1/D_1} \mathbf{s}) \xi(d(\mathbf{s},\mathbf{t})),$$

where the function $K$ is our Kernel, which fulfills the following condition:

CONDITION 5 (Shape of $K$). The kernel $K : \mathbb{R}^{D_1} \to \mathbb{R}_+$ satisfies:

(i) $K(\mathbf{s}) = 0$ for $\mathbf{s} \notin [-1,1)^{D_1}$;
(ii) $K|_{[-1,1)^{D_1}}$ is Lipschitz (w.r.t. $d_0$ restricted to $\mathbb{R}^{D_1}$) with constant $l(K)$;
(iii) $\int K(\mathbf{s})\,d\mathbf{s} = 1$;
(iv) $\int K(\mathbf{s})\mathbf{s}\,d\mathbf{s} = \mathbf{0}$.

Note that $K$ does not have to be continuous on the boundary of $[-1,1)^{D_1}$, and that it is reasonable to choose a Kernel $K$ that is radially symmetric (or at least an even function in each coordinate), in which case Condition 5(iv) is satisfied. We now write

$$f(\mathbf{x}) := 2^{D_1} K(\mathbf{s}) \cdot \mathbb{1}_{[-1,1)^{D_2}}(\mathbf{t}) \qquad \text{for } \mathbf{x} := (\mathbf{s},\mathbf{t}) \in \mathbb{R}^{D_1} \times \mathbb{R}^{D_2} = \mathbb{R}^D,$$

so that $f|_J$ is Lipschitz (w.r.t. $d_0$ on $\mathbb{R}^D$) with constant $2^{D_1} l(K)$; by the transformation theorem for integrals, we obtain

$$\hat{p}_\xi(\mathbf{0}) = \frac{1}{|J_T|} \int_{\mathbb{R}^D} f(\mathbf{x}) \xi \theta_T^{-1}(d\mathbf{x}).$$

The way is now clear for the application of Theorem 2.A. Our primary goal will be to estimate a probability distance $d$ between the distribution of our estimator $\hat{p}_\xi(\mathbf{0})$ and the distribution that is concentrated at the true value $p(\mathbf{0})$. To do this, we will first estimate $d(\mathcal{L}(\hat{p}_\xi(\mathbf{0})), \mathcal{L}(\hat{p}_\eta(\mathbf{0})))$ with the aid of Theorem 2.A, and then utilize the excellent properties of Poisson point processes to obtain an upper bound for $d(\mathcal{L}(\hat{p}_\eta(\mathbf{0})), \delta_{p(\mathbf{0})})$. The two

BOUNDS FOR POINT PROCESS APPROXIMATIONS 29

corresponding results are contained in the following theorems. For the distance $d$, we choose the bounded Wasserstein distance, as defined in Section 1, because the other distances that we have used so far are too strong to be useful: $d_{\mathrm{TV}}(\mathcal{L}(\hat{p}_\xi(\mathbf{0})), \delta_{p(\mathbf{0})})$ is generally too big, and is even always equal to 1 whenever $\hat{p}_\xi(\mathbf{0})$ is a continuous random variable, because then

$$1 \geq d_{\mathrm{TV}}(\mathcal{L}(\hat{p}_\xi(\mathbf{0})), \delta_{p(\mathbf{0})}) \geq |\mathbb{P}[\hat{p}_\xi(\mathbf{0}) = p(\mathbf{0})] - \mathbb{P}[p(\mathbf{0}) = p(\mathbf{0})]| = 1;$$

and for the Wasserstein distance $d_{\mathbf{w}}(\mathcal{L}(\hat{p}_\xi(\mathbf{0})), \mathcal{L}(\hat{p}_\eta(\mathbf{0})))$, there seem to be unsurmountable difficulties in obtaining a useful upper bound in Theorem 3.A.

THEOREM 3.A. *Suppose that the prerequisites of Section 1 hold, including the Conditions 1, 2, 3$\rho$, as well as the additional Conditions 4′ and 5. Let $\iota > 0$, and for $T \geq 1$, let $m := m(T) \in \mathbb{Z}_+$, $h(T) \geq 1$ and also $w(T) = O(T^{\delta^*})$ for $T \to \infty$ with $\delta^* \in (0,1)$. Then*

$$d_{\mathrm{BW}}(\mathcal{L}(\hat{p}_\xi(\mathbf{0})), \mathcal{L}(\hat{p}_\eta(\mathbf{0})))$$
$$\leq \left(\frac{l(K)}{2^{D_2}} \frac{w(T)}{T} M + 1\right) d_2(\mathcal{L}(\xi\theta_T^{-1}|_J), \mathcal{L}(\eta\theta_T^{-1}|_J)) + 2^{D_1} l(K) \delta_T(M)$$
$$= O\left(\frac{1}{h(T)^{1/D_1}}, \frac{1}{T^{1/D_2}}, \log^\uparrow\left(\frac{T}{w(T)}\right) \frac{m^{D_2}+1}{w(T)}, \frac{T}{w(T)} \check{\alpha}\left(\frac{2^{D_2}}{w(T)h(T)}\right),\right.$$
$$\left. \log^\uparrow\left(\frac{T}{w(T)}\right) \check{\alpha}\left(\frac{2^D(2m+1)^{D_2}}{w(T)}\right), \sqrt{Th(T)}\check{\beta}(m)\right)$$
$$\text{for } T \to \infty,$$

*where $M := M(T) \in \mathbb{N}^*$ with $M \geq 3\nu(J_T)$ arbitrary and*

$$\delta_T(M) = 2\kappa \frac{\nu(J_T)^M}{M!} e^{-\nu(J_T)},$$

*which decays exponentially in $M$ as $T$ tends to infinity. Thus, we obtain the same order for the upper bound as in Theorem 2.A*

REMARK 3.B. The upper bound given in Theorem 3.A remains true for general $w(T) = O(T)$. However, if $w(T)$ goes to infinity at a rate that is too close to $T$, then $M(T)$ has to be chosen to grow somewhat faster than $T/w(T)$, and then the order of the upper bound is a little worse (by a logarithmic factor in $T$) than the one stated in Theorem 3.A.

PROOF OF THEOREM 3.A. Let $\xi' \sim \mathcal{L}(\xi)$, $\eta' \sim \mathcal{L}(\eta) = \mathrm{Po}(\nu)$, and $X := \hat{p}_{\xi'}(\mathbf{0})$, $Y := \hat{p}_{\eta'}(\mathbf{0})$. Then we have

$$d_{\mathrm{BW}}(\mathcal{L}(\hat{p}_\xi(\mathbf{0})), \mathcal{L}(\hat{p}_\eta(\mathbf{0}))) = \sup_{g \in \mathcal{F}_{\mathrm{BW}}} |\mathbb{E}g(X) - \mathbb{E}g(Y)|$$

with

$$
\begin{aligned}
|\mathbb{E}g(X) - \mathbb{E}g(Y)| \\
&\leq \mathbb{E}(|g(X) - g(Y)|\mathbb{1}_{\{\xi'\theta_T^{-1}(J)=\eta'\theta_T^{-1}(J)\}}) \\
&\quad + \mathbb{E}(|g(X) - g(Y)|\mathbb{1}_{\{\xi'\theta_T^{-1}(J)\neq\eta'\theta_T^{-1}(J)\}}) \\
&\leq \mathbb{E}(|X - Y|\mathbb{1}_{\{\xi'\theta_T^{-1}(J)=\eta'\theta_T^{-1}(J)\}}) + \mathbb{P}[\xi'\theta_T^{-1}(J) \neq \eta'\theta_T^{-1}(J)]
\end{aligned}
\tag{3.1}
$$

for every $g$ in $\mathcal{F}_{\mathrm{BW}}$. For the first summand, we obtain

$$
\begin{aligned}
&\mathbb{E}(|X - Y|\mathbb{1}_{\{\xi'\theta_T^{-1}(J)=\eta'\theta_T^{-1}(J)\}}) \\
&= \mathbb{E}\bigg(\frac{1}{|J_T|}\bigg|\int_{\mathbb{R}^D} f(\mathbf{x})\xi'\theta_T^{-1}(d\mathbf{x}) - \int_{\mathbb{R}^D} f(\mathbf{x})\eta'\theta_T^{-1}(d\mathbf{x})\bigg|\mathbb{1}_{\{\xi'\theta_T^{-1}(J)=\eta'\theta_T^{-1}(J)\}}\bigg) \\
&\leq 2^{D_1}l(K)\mathbb{E}\bigg(\frac{\eta'\theta_T^{-1}(J)}{|J_T|}d_1(\xi'\theta_T^{-1}|_J, \eta'\theta_T^{-1}|_J)\bigg),
\end{aligned}
$$

the latter inequality by the definition of the $d_1$-distance and because $f|_J$ is Lipschitz. Next we utilize the fact that since $\eta'\theta_T^{-1}(J)$ is Poisson distributed with parameter $\nu_T := \nu(J_T)$, it exceeds a certain bound $M := M(T) \in \mathbb{N}^*$ with $M + 1 \geq 2\nu_T$ only with very small probability. As noted in Barbour, Holst and Janson (1992), Proposition A.2.3, the relation

$$
\mathbb{P}[\mathrm{Po}(\nu_T) \geq M] \leq \frac{M+1}{M+1-\nu_T}\mathbb{P}[\mathrm{Po}(\nu_T) = M] \leq 2\frac{\nu_T^M}{M!}e^{-\nu_T}
$$

holds, and, thus,

$$
\begin{aligned}
&\mathbb{E}\bigg(\frac{\eta'\theta_T^{-1}(J)}{|J_T|}d_1(\xi'\theta_T^{-1}|_J, \eta'\theta_T^{-1}|_J)\bigg) \\
&\leq \mathbb{E}\bigg(\frac{M}{|J_T|}d_1(\xi'\theta_T^{-1}|_J, \eta'\theta_T^{-1}|_J)\mathbb{1}_{\{\eta'\theta_T^{-1}(J)\leq M\}}\bigg) \\
&\quad + \mathbb{E}\bigg(\frac{\eta'\theta_T^{-1}(J)}{|J_T|}\mathbb{1}_{\{\eta'\theta_T^{-1}(J)>M\}}\bigg) \\
&\leq \frac{M}{|J_T|}\mathbb{E}(d_1(\xi'\theta_T^{-1}|_J, \eta'\theta_T^{-1}|_J)) + \frac{\nu_T}{|J_T|}\mathbb{P}[\eta'\theta_T^{-1}(J) \geq M] \\
&\leq \frac{1}{2^D}\frac{w(T)}{T}M\mathbb{E}(d_1(\xi'\theta_T^{-1}|_J, \eta'\theta_T^{-1}|_J)) + \delta_T(M),
\end{aligned}
$$

where we use the notation

$$
\delta_T(M) = 2\kappa\frac{\nu_T^M}{M!}e^{-\nu_T}.
$$

Furthermore, for $M \geq 3\nu_T$, the DeMoivre–Stirling formula gives

$$
\delta_T(M) \leq \mathrm{const}\cdot\bigg(\frac{\nu_T}{M}\bigg)^M e^{M-\nu_T} \leq \mathrm{const}\cdot\bigg(\frac{e}{3}\bigg)^M e^{-\nu_T}.
$$

The second summand from (3.1) is estimated as

$$\mathbb{P}[\xi'\theta_T^{-1}(J) \neq \eta'\theta_T^{-1}(J)] = \mathbb{E}[d_1(\xi'\theta_T^{-1}|_J, \eta'\theta_T^{-1}|_J)\mathbb{1}_{\{\xi'\theta_T^{-1}(J) \neq \eta'\theta_T^{-1}(J)\}}]$$

$$\leq \mathbb{E}d_1(\xi'\theta_T^{-1}|_J, \eta'\theta_T^{-1}|_J).$$

Hence, we obtain altogether in (3.1),

$$|\mathbb{E}g(X) - \mathbb{E}g(Y)|$$
$$\leq \left(\frac{l(K)}{2^{D_2}}\frac{w(T)}{T}M + 1\right)\mathbb{E}(d_1(\xi'\theta_T^{-1}|_J, \eta'\theta_T^{-1}|_J)) + 2^{D_1}l(K)\delta_T(M)$$

for every $g \in \mathcal{F}_{\mathrm{BW}}$ and every pair of random variables $\xi'$, $\eta'$ with $\xi' \sim \mathcal{L}(\xi)$, $\eta' \sim \mathcal{L}(\eta)$. Forming the infimum over $\xi'$ and $\eta'$ yields on the right-hand side the $d_2$-distance ($\theta_T$ is bijective), and forming the supremum over $g$ on the left-hand side, the bounded Wasserstein distance. Thus, we obtain the statement. □

The second result that was discussed above is contained in the next theorem. We write $\|\cdot\|_2$ for the $L_2$-norm with respect to the Lebesgue measure on $\mathbb{R}^{D_1}$.

THEOREM 3.C. *Suppose that the prerequisites of Section* 1 *hold, including the Conditions* 1, 2, 3$\rho$, *as well as the additional Conditions* 4′ *and* 5. *Let* $\iota > 0$, *and for* $T \geq 1$, *let* $m := m(T) \in \mathbb{Z}_+$, $h(T) \geq 1$ *and also* $w(T) = O(T^{\delta^*})$ *for* $T \to \infty$ *with* $\delta^* \in (0,1)$. *Then*

$$d_{\mathrm{BW}}(\mathcal{L}(\hat{p}_\xi(\mathbf{0})), \delta_{p(\mathbf{0})})$$
$$\leq d_{\mathrm{BW}}(\mathcal{L}(\hat{p}_\xi(\mathbf{0})), \mathcal{L}(\hat{p}_\eta(\mathbf{0})))$$
$$+ \sqrt{\frac{\kappa}{2^{D_2}}}\|K\|_2\sqrt{\frac{w(T)}{T}} + \frac{L'}{w(T)^{2/D_1}} + o\left(\frac{1}{w(T)^{2/D_1}}\right)$$
$$= O\left(\sqrt{\frac{w(T)}{T}}, \frac{1}{w(T)^{2/D_1}}, \frac{1}{h(T)^{1/D_1}}, \frac{1}{T^{1/D_2}},\right.$$
$$\log^\uparrow\left(\frac{T}{w(T)}\right)\frac{m^{D_2}+1}{w(T)}, \frac{T}{w(T)}\breve{\alpha}\left(\frac{2^{D_2}}{w(T)h(T)}\right),$$
$$\left.\log^\uparrow\left(\frac{T}{w(T)}\right)\breve{\alpha}\left(\frac{2^D(2m+1)^{D_2}}{w(T)}\right), \sqrt{Th(T)}\breve{\beta}(m)\right)$$

for $T \to \infty$,

*where* $L'$ *is a nonnegative constant (depending on* $p$ *and* $K$); *if* $K$ *possesses certain symmetry properties (especially if* $K$ *is radially symmetric), we can*

*write*

$$L' := \tfrac{1}{2}\Delta p(\mathbf{0}) \int s_1^2 K(\mathbf{s})\lambda^{D_1}(d\mathbf{s}),$$

*where $\Delta$ denotes the $D_1$-dimensional Laplace operator.*

PROOF. Due to Theorem 3.A we only have to estimate $d_{\mathrm{BW}}(\mathcal{L}(\hat{p}_\eta(\mathbf{0})), \delta_{p(\mathbf{0})})$ for $\eta \sim \mathrm{Po}(\nu)$. We decompose this distance as

$$\begin{aligned}
d_{\mathrm{BW}}(\mathcal{L}(\hat{p}_\eta(\mathbf{0})), \delta_{p(\mathbf{0})}) &\leq d_{\mathrm{BW}}(\mathcal{L}(\hat{p}_\eta(\mathbf{0})), \delta_{\mathbb{E}\hat{p}_\eta(\mathbf{0})}) + d_{\mathrm{BW}}(\delta_{\mathbb{E}\hat{p}_\eta(\mathbf{0})}, \delta_{p(\mathbf{0})}) \\
&\leq \mathbb{E}|\hat{p}_\eta(\mathbf{0}) - \mathbb{E}\hat{p}_\eta(\mathbf{0})| + |\mathbb{E}\hat{p}_\eta(\mathbf{0}) - p(\mathbf{0})| \\
&\leq \mathrm{sd}(\hat{p}_\eta(\mathbf{0})) + \mathrm{bias}(\hat{p}_\eta(\mathbf{0})).
\end{aligned}$$

For the standard deviation we obtain

$$\begin{aligned}
\mathrm{sd}(\hat{p}_\eta(\mathbf{0})) &= \sqrt{\mathrm{var}\left(\frac{1}{|J_T|}\int_{\mathbb{R}^D} f(\mathbf{x})\eta\theta_T^{-1}(d\mathbf{x})\right)} \\
&= \frac{1}{|J_T|}\sqrt{\int_{\mathbb{R}^D} f^2(\mathbf{x})\nu\theta_T^{-1}(d\mathbf{x})} \\
&\leq \frac{1}{|J_T|}\sqrt{\kappa_T\left(\frac{1}{w(T)}\int_{\mathbb{R}^{D_1}} 2^{2D_1}K^2(\mathbf{s})\lambda^{D_1}(d\mathbf{s})\right)\mu_2([-T^{1/D_2}, T^{1/D_2})^{D_2})} \\
&\leq \sqrt{\frac{\kappa}{2^{D_2}}}\|K\|_2\sqrt{\frac{w(T)}{T}},
\end{aligned}$$

where the second and third steps are applications of Campbell's theorem for the variance of an integral w.r.t. a Poisson point process [see Kingman (1993)] and Fubini's theorem, respectively [note that $(\lambda^{D_1} \otimes \mu_2)\theta_T^{-1} = \frac{1}{w(T)}\lambda^{D_1} \otimes \mu_2(T^{1/D_2}I_{D_2})$, where $I_{D_2}\colon\mathbb{R}^{D_2} \to \mathbb{R}^{D_2}$ is the identity]. An application of Campbell's theorem for the expectation [see Kingman (1993)] and Fubini's theorem again then yields

$$\begin{aligned}
\mathbb{E}\hat{p}_\eta(\mathbf{0}) &= \frac{1}{|J_T|}\int_{\mathbb{R}^D} f(\mathbf{x})\nu\theta_T^{-1}(d\mathbf{x}) \\
&= \frac{1}{|J_T|}\left(\frac{1}{w(T)}\int_{\mathbb{R}^{D_1}} 2^{D_1}K(\mathbf{s})p\left(\frac{1}{w(T)^{1/D_1}}\mathbf{s}\right)\lambda^{D_1}(d\mathbf{s})\right) \\
&\quad \times \mu_2([-T^{1/D_2}, T^{1/D_2})^{D_2}) \\
&= \int_{\mathbb{R}^{D_1}} K(\mathbf{s})p\left(\frac{1}{w(T)^{1/D_1}}\mathbf{s}\right)\lambda^{D_1}(d\mathbf{s}).
\end{aligned}$$

BOUNDS FOR POINT PROCESS APPROXIMATIONS 33Thus, we obtain for the bias

$$|\mathbb{E}\hat{p}_\eta(\mathbf{0}) - p(\mathbf{0})|$$
$$= \left|\int_{[-1,1)^{D_1}} K(\mathbf{s})\left(p\left(\frac{1}{w(T)^{1/D_1}}\mathbf{s}\right) - p(\mathbf{0})\right)\lambda^{D_1}(d\mathbf{s})\right|$$
$$\leq \left|\int_{[-1,1)^{D_1}} K(\mathbf{s})\frac{1}{w(T)^{1/D_1}}\partial p(\mathbf{0})\mathbf{s}\,\lambda^{D_1}(d\mathbf{s})\right|$$
$$+ \left|\int_{[-1,1)^{D_1}} K(\mathbf{s})\frac{1}{2w(T)^{2/D_1}}\partial^2 p(\mathbf{0})(\mathbf{s},\mathbf{s})\lambda^{D_1}(d\mathbf{s})\right|$$
$$+ \int_{[-1,1)^{D_1}} K(\mathbf{s})\frac{1}{2w(T)^{2/D_1}}$$
$$\times \max_{0\leq h\leq 1}\left\|\partial^2 p\left(h\frac{1}{w(T)^{1/D_1}}\mathbf{s}\right) - \partial^2 p(\mathbf{0})\right\||\mathbf{s}|^2\lambda^{D_1}(d\mathbf{s})$$

by Taylor's approximation, where $\|\cdot\|$ is the standard norm for bilinear forms on $\mathbb{R}^{D_1}$. Of the last three summands, the first is always zero because of Condition 5(iv), the second can be estimated by $L'\frac{1}{w(T)^{2/D_1}}$ with a constant $L'$, which for "nice" Kernels (e.g., if $K$ is radially symmetric) can be written as

$$L' = \tfrac{1}{2}\Delta p(\mathbf{0})\int s_1^2 K(\mathbf{s})\lambda^{D_1}(d\mathbf{s}),$$

and the third is of order $o(\frac{1}{w(T)^{2/D_1}})$ because of the continuity of $\partial^2 p$ at $\mathbf{0}$. Thus,

$$\text{bias}(\hat{p}_\eta(\mathbf{0})) \leq L'\frac{1}{w(T)^{2/D_1}} + o\left(\frac{1}{w(T)^{2/D_1}}\right). \qquad \square$$

Once more we formulate the conditions under which the upper bound goes to zero.

COROLLARY 3.D (Convergence to zero in Theorem 3.C). *Suppose that the prerequisites of Theorem 3.C hold. Furthermore, suppose that $w(T) \geq kT^\delta$ for $k > 0$, $\delta \in (0,1)$ and that*

$$\check{\alpha}(v) = O(v^r) \qquad \text{for } v \to 0 \text{ with } r > 0,$$
$$\check{\beta}(u) = O\left(\frac{1}{u^{(1+s)D_2/2}}\right) \qquad \text{for } u \to \infty \text{ with } 1+s > \max\left(\frac{1-\delta}{\delta}\frac{1+r}{r}, \frac{1}{\delta}\right).$$

*Then*

$$d_{\text{BW}}(\mathcal{L}(\hat{p}_\xi(\mathbf{0})), \delta_{p(\mathbf{0})}) \to 0 \qquad \text{for } T \to \infty,$$

and, therefore, since the $d_{\mathrm{BW}}$-distance metrizes convergence in distribution [see Dudley (1989), Theorem 11.3.3] and since $\delta_{p(\mathbf{0})}$ is the distribution of a constant, we obtain

$$\hat{p}_\xi(\mathbf{0}) \xrightarrow{P} p(\mathbf{0}) \qquad \text{for } T \to \infty,$$

that is, the consistency of the estimator $\hat{p}_\xi(\mathbf{0})$.

REMARK 3.E. The consistency of $\hat{p}_\xi(\mathbf{0})$ was already obtained as a consequence of Theorem 2.5 in Ellis (1991) under conditions that were similar, but for the most part somewhat more general. So Corollary 3.D is not so much a new result, but rather a crosscheck on the suitability of the explicit upper bound obtained in Corollary 3.C.

PROOF. Let $M := \lceil 3\nu(J_T) \rceil$ in Theorem 3.A. We then get immediately by applying Theorems 3.C and 3.A and Corollary 2.B that $d_{\mathrm{BW}}(\mathcal{L}(\hat{p}_\xi(\mathbf{0})), \delta_{p(\mathbf{0})})$ converges to zero. □

3.2. *Testing for long range dependence.* Suppose $\xi$ is a stationary point process on $\mathbb{R}^D$ with expectation measure $\nu = \ell \cdot \lambda^D$ ($\ell$ known or estimated) which satisfies the conditions of Section 1, except for Condition 3. We would like to test from a single realization of $\xi$ if there is important long range dependence in the $\mathbb{R}^{D_2}$-directions or not (our null hypothesis). "No important long range dependence" means here that Condition $3x$ is satisfied for given $x \in \{\beta, \rho, \varphi\}$ and $\breve{\beta}$, corresponding to the minimal mixing rate one wants to test for. For the sake of illustration, think of the $\mathbb{R}^{D_1}$-direction(s) as time and the $\mathbb{R}^{D_2}$-directions as space. Imagine that for fixed $T \geq 1$, the points of $\xi$ in $J_T$ denote the times and locations of incidences of a certain rare disease, which is observed in a large area (e.g., a country or a continent) over a relatively short period of time (e.g., some months or a year).

Under the null hypothesis, by Theorem 2.A, respectively, Theorem 2.D, the distribution of $\xi\theta_T^{-1}|_J$ will be close to the distribution of $\eta\theta_T^{-1}|_J$, which here is just the homogeneous Poisson process on $J$ with intensity $(T/w(T)) \cdot \ell$. There are various reasonable statistics for testing the hypothesis of "complete spatial randomness" in point patterns; one such statistic, $U : \mathcal{M}_p \to \mathbb{R}$, is the average nearest neighbor distance in the data, which can be shown to be Lipschitz continuous with respect to the $d_1$-distance with a Lipschitz constant that we denote by $L_D$.

We wish to find an approximate critical value $t_\alpha$ for, say, a one-sided test of size $\alpha$ of the null hypothesis against an aggregated alternative (i.e., the alternative that there is a certain amount of "long range" clustering), using the statistic $\tilde{U}$, where $\tilde{U}(\rho) := U(\rho\theta_T^{-1}|_J)$ for every point measure $\rho$ on $\mathbb{R}^D$. To do so, fix $K > 0$ and choose $t_\alpha$ so that

$$\mathbb{E} f_{t_\alpha, K}(\tilde{U}(\eta)) + KL_D \cdot \varepsilon = \alpha,$$

where $\varepsilon$ is our upper bound for $d_2(\mathcal{L}(\xi\theta_T^{-1}|_J), \mathcal{L}(\eta\theta_T^{-1}|_J))$, and

$$f_{t,K}(x) := \begin{cases} 1, & \text{if } x \leq t, \\ 1 - K(x-t), & \text{if } t \leq x \leq t + \frac{1}{K}, \\ 0, & \text{if } x \geq t + \frac{1}{K}, \end{cases}$$

is a $K$-Lipschitz approximation of the indicator $\mathbb{1}_{(-\infty,t]}$. This yields

$$0 \leq \alpha - \mathbb{P}[\tilde{U}(\xi) < t_\alpha]$$
$$\leq \mathbb{E}f_{t_\alpha,K}(\tilde{U}(\eta)) - \mathbb{E}f_{(t_\alpha - 1/K),K}(\tilde{U}(\eta)) + 2KL_D \cdot \varepsilon.$$

Thus, if $\varepsilon$ is very small (i.e., the conditions for Theorem 2.A, resp. Theorem 2.D, are strong enough), a large $K$ can be chosen, and, consequently, we can adjust the size of our test to be only slightly below $\alpha$.

It should be noted that the distribution of $\tilde{U}(\eta)$ is not known, but it can be simulated very easily. Also, there are good normal approximations of $\mathcal{L}(\tilde{U}(\eta)||\eta| = N)$ for $N$ not too small which can be of use. See Ripley [(1981), Section 8.2] for further details.

## APPENDIX: LOCAL STEIN THEOREMS

The central results of this article are achieved by applying estimates that were obtained in one or another form by Stein's method. Since it is far beyond the scope of this article to summarize in detail the classical Stein–Chen method (Stein's method for the approximation of a sum of indicator random variables by a Poisson random variable) or what in this article is sometimes called the "generalized Stein–Chen method" (Stein's method for the approximation of an indicator point process by a discrete Poisson point process), we only present very briefly the required results. The proofs of these results and the method behind them, as well as a wealth of related material, can be found in Barbour, Holst and Janson (1992).

Let $\Gamma$ be any finite nonempty index set and $(I_i)_{i \in \Gamma}$ a sequence of indicator random variables with a local dependence property, that is, for every $i \in \Gamma$, the set $\Gamma_i := \Gamma \setminus \{i\}$ can be partitioned as $\Gamma_i = \Gamma_i^s \dot{\cup} \Gamma_i^w$ into a set $\Gamma_i^s$ of indices $j$, for which $I_j$ depends "strongly" on $I_i$, and a set $\Gamma_i^w$ of indices $j$, for which $I_j$ depends "weakly" on $I_i$. Herein, the terms "strongly" and "weakly" are not meant as a restriction to the partition of $\Gamma_i$, but serve only illustrative purposes. The same holds true for the term "local dependence," which does not have to possess any representation in the spatial structure of $\Gamma$ (in our applications in Section 2 it always does, though). We now write $Z_i := \sum_{j \in \Gamma_i^s} I_j$, $Y_i := \sum_{j \in \Gamma_i^w} I_j$, $p_i := \mathbb{E}I_i > 0$ (w.l.o.g.) for every $i \in \Gamma$ and set $W := \sum_{i \in \Gamma} I_i$, $\lambda := \mathbb{E}W = \sum_{i \in \Gamma} p_i$. Furthermore, we choose arbitrary points $(\alpha_i)_{i \in \Gamma}$ in any desired complete, separable metric space $(\mathcal{X}, d_0)$ with $d_0 \leq 1$ and set $\Xi := \sum_{i \in \Gamma} I_i \delta_{\alpha_i}$.

**A.1. Poisson approximation of the distribution of the sum $W$ of indicators.** By applying the classical Stein–Chen method [see Chen (1975)] the following result is obtained.

THEOREM A.A (Local Stein–Chen theorem for sums of indicators). *With the above definitions, we have*

$$d_{\mathrm{TV}}(\mathcal{L}(W), \mathrm{Po}(\lambda))$$
$$\leq \min\left(1, \frac{1}{\lambda}\right) \sum_{i \in \Gamma} (p_i^2 + p_i \mathbb{E} Z_i + \mathbb{E}(I_i Z_i)) + \min\left(1, \frac{1}{\sqrt{\lambda}}\right) \sum_{i \in \Gamma} e_i,$$

*where*

$$e_i = \mathbb{E}|\mathbb{E}(I_i | (I_j : j \in \Gamma_i^w)) - p_i| = 2 \max_{B \in \sigma(I_j : j \in \Gamma_i^w)} |\mathrm{cov}(I_i, \mathbb{1}_B)|.$$

PROOF. See Barbour, Holst and Janson (1992), Theorem 1.A. □

REMARK A.B. The order of the upper bound in Theorem A.A cannot generally be improved. See Barbour, Holst and Janson (1992), Chapter 3.

The Stein–Chen method is by no means restricted to approximating sums of indicator random variables. For instance, as far as $\mathbb{Z}_+$-valued random variables are concerned, one might also consider the case where $W$ is itself Poisson distributed with some parameter $\mu > 0$.

PROPOSITION A.C. *Let $\lambda, \mu > 0$. Then*

$$d_{\mathrm{TV}}(\mathrm{Po}(\lambda), \mathrm{Po}(\mu)) \leq \min\left(1, \frac{1}{\sqrt{\lambda}}, \frac{1}{\sqrt{\mu}}\right) \cdot |\lambda - \mu|.$$

PROOF. This proposition is a special case of Barbour, Holst and Janson (1992), Theorem 1.C(i). However, the result can be obtained very easily by direct calculation, using the Stein–Chen method. □

**A.2. Poisson process approximation of the distribution of the indicator point process $\Xi$.** By applying a natural generalization of the Stein–Chen method as in Barbour and Brown (1992), the following result is obtained.

THEOREM A.D (Local Stein theorem for indicator point processes). *With the above definitions and $\boldsymbol{\pi} := \sum_{i \in \Gamma} p_i \delta_{\alpha_i}$, we have*

$$d_2(\mathcal{L}(\Xi), \mathrm{Po}(\boldsymbol{\pi}))$$
$$\leq \left\{ 1 \wedge \frac{2}{\lambda}\left(1 + 2\log^+\left(\frac{\lambda}{2}\right)\right) \right\} \sum_{i \in \Gamma} (p_i^2 + p_i \mathbb{E} Z_i + \mathbb{E}(I_i Z_i))$$
$$+ \left(1 \wedge 1.65 \frac{1}{\sqrt{\lambda}}\right) \sum_{i \in \Gamma} e_i,$$

*where*

$$e_i = \mathbb{E}|\mathbb{E}(I_i|(I_j; j \in \Gamma_i^w)) - p_i| = 2 \max_{B \in \sigma(I_j; j \in \Gamma_i^w)} |\operatorname{cov}(I_i, \mathbb{1}_B)|.$$

PROOF. See Barbour, Holst and Janson (1992), Theorem 10.F. □

REMARK A.E. Note that the upper bound in Theorem A.D depends neither on the points $\alpha_i$, $i \in \Gamma$, nor on the specific choice of the metric $d_0$, as long as it is bounded by 1.

**Acknowledgment.** The author wishes to thank Andrew Barbour for contributing many helpful comments and suggestions to this article.

## REFERENCES

BARBOUR, A. D. and BROWN, T. C. (1992). Stein's method and point process approximation. *Stochastic Process. Appl.* **43** 9–31. MR1190904

BARBOUR, A. D., HOLST, L. and JANSON, S. (1992). *Poisson Approximation.* Oxford Univ. Press. MR1163825

BOROVKOV, A. A. (1996). Asymptotic expansions for functionals of dilation of point processes. *J. Appl. Probab.* **33** 573–591. MR1385366

CHEN, L. H. Y. (1975). Poisson approximation for dependent trials. *Ann. Probab.* **3** 534–545. MR428387

DALEY, D. J. (1974). Various concepts of orderliness for point processes. In *Stochastic Geometry* (E. F. Harding and D. G. Kendall, eds.) 148–161. Wiley, New York. MR380976

DOUKHAN, P. (1994). *Mixing*: *Properties and Examples.* Springer, New York. MR1312160

DUDLEY, R. M. (1989). *Real Analysis and Probability.* Wadsworth & Brooks/Cole, Pacific Grove, CA. MR982264

ELLIS, S. P. (1986). A limit theorem for spatial point processes. *Adv. in Appl. Probab.* **18** 646–659. MR857323

ELLIS, S. P. (1991). Density estimation for point processes. *Stochastic Process. Appl.* **39** 345–358. MR1136254

KALLENBERG, O. (1986). *Random Measures*, 4th ed. Academic Press, New York. MR854102

KINGMAN, J. F. C. (1993). *Poisson Processes.* Oxford Univ. Press. MR1207584

RACHEV, S. T. (1984). The Monge–Kantorovich mass transference problem and its stochastic applications. *Theory Probab. Appl.* **29** 647–676. MR773434

RIPLEY, B. D. (1981). *Spatial Statistics.* Wiley, New York. MR624436

SILVERMAN, B. W. (1986). *Density Estimation for Statistics and Data Analysis.* Chapman and Hall, London. MR848134

INSTITUTE OF MATHEMATICS
UNIVERSITY OF ZÜRICH
WINTERTHURERSTRASSE 190
CH-8057 ZÜRICH
SWITZERLAND
E-MAIL: schumi@amath.unizh.ch